\documentclass[final]{siamltex}
\usepackage[utf8]{inputenc}
\usepackage{amssymb}
\usepackage{amsfonts}
\usepackage{mathrsfs}
\usepackage{amsmath}
\usepackage{graphicx}
\usepackage{dsfont}
\usepackage{float}
\usepackage{diagbox} 
\newcommand{\norm}[1]{\left\lVert#1\right\rVert}
\usepackage{graphics}
\usepackage[all]{xy}
\usepackage[margin=1.4in,marginparwidth=1in,marginparsep=.05in]{geometry}
\PassOptionsToPackage{hyphens}{url}
\usepackage{tikz-cd}
\usepackage{braket}
\usepackage{xcolor}
\usepackage{verbatim}
\usepackage{float}
\usepackage{mathrsfs}
\usepackage{extarrows}
\usepackage{tikz}
\usetikzlibrary{shapes,arrows,positioning}
\usepackage{mathrsfs}
\usepackage{mathtools}

\usepackage{subfig}
\usepackage{hyperref}
\usepackage{wrapfig}

\newtheorem{remark}[theorem]{Remark}

\numberwithin{equation}{section}

\let\pa\partial

\newcommand{\R}{\mathbb R}
\newcommand{\bA}{\mathbf A}

\newcommand{\bD}{\mathbf D}

\newcommand{\bI}{\mathbf I}
\newcommand{\bJ}{\mathbf J}
\newcommand{\bP}{\mathbf P}
\newcommand{\bG}{\mathbf G}

\newcommand{\bV}{\mathbf V}

\newcommand{\ba}{\mathbf a}
\newcommand{\bb}{\mathbf b}
\newcommand{\bc}{\mathbf c}
\newcommand{\bg}{\mathbf g}

\newcommand{\blf}{\mathbf f}
\newcommand{\bn}{\mathbf n}
\newcommand{\be}{\mathbf e}

\newcommand{\bu}{\mathbf u}
\newcommand{\bU}{\mathbf U}
\newcommand{\bv}{\mathbf v}
\newcommand{\bw}{\mathbf w}

\newcommand{\bx}{\mathbf x}

\newcommand{\bbf}{\mathbf f}

\newcommand{\T}{\mathcal T}

\newcommand{\Div}{\mathop{\rm div}}
\newcommand{\divG}{{\mathop{\,\rm div}}_{\Gamma}}
\newcommand{\gradG}{\nabla_{\Gamma}}
\newcommand{\nablaG}{\nabla_{\Gamma}}

\newcommand{\OGamma}{\Omega^\Gamma_h}

\renewcommand{\div}{\textrm{div}\ \!}

\newcommand{\tr}{{\rm tr}}

\DeclareGraphicsExtensions{.pdf,.eps,.ps,.eps.gz,.ps.gz,.eps.Y}

\newcommand{\la}{\left\langle}
\newcommand{\ra}{\right\rangle}

\newcommand{\bsigma}{\boldsymbol{\sigma}}

\def\cl {\nonumber \\}
\def\el {\nonumber }

\newcommand{\compositenorm}[1]{{\left\vert\kern-0.25ex\left\vert\kern-0.25ex\left\vert #1
    \right\vert\kern-0.25ex\right\vert\kern-0.25ex\right\vert}}

\usepackage{accents}
\renewcommand*{\dot}[1]{%
   \accentset{\mbox{\large\bfseries .}}{#1}}

\title{A finite element method for two-phase flow with material viscous interface}
\author{Maxim Olshanskii$^{*}$, Annalisa Quaini$^{*}$ and Qi Sun\thanks{Department of Mathematics, University of Houston, 3551 Cullen Blvd, Houston TX 77204 (\texttt{maolshanskiy@uh.edu; aquaini@uh.edu; qsun5@uh.edu})\\
This work has been partially supported by US National Science Foundation (NSF) through grants DMS-1620384, DMS-1953535, DMS-2011444.
}
}
\date{\today}

\begin{document}

\maketitle


\begin{abstract}
This paper studies a model of two-phase flow with an immersed material viscous interface and a finite element method for  numerical solution of the resulting system of PDEs. The interaction between the bulk and surface media is characterized by no-penetration and slip with friction interface conditions.  
The system is shown to be dissipative and a model stationary problem is proved to be well-posed.
The  finite element method applied in this paper belongs to a family of unfitted discretizations. The 
performance of the method when model and discretization parameters vary is assessed. Moreover,
an iterative procedure based on the splitting of the system into bulk and surface problems is introduced and studied numerically.
\end{abstract}

\begin{keywords}
	two-phase flow, surface Stokes equations, finite elements, cutFEM, TraceFEM
\end{keywords}

\section{Introduction}\label{sec:intro}
We consider two immiscible, viscous, and incompressible fluids 
separated by a viscous  inextensible material interface modeled as a Boussinesq--Scri\-ven surface fluid.
The following coupling conditions are prescribed 
between the bulk two-phase flow and the surface fluid: (i) the immiscibility condition, i.e.~the bulk fluid does not penetrate through
the interface; (ii) slip with friction between the bulk fluid and the viscous interface; and (iii) the
load exerted from the bulk fluid onto the surface fluid defined by the jump of the normal stress across the interface.
To the best of our knowledge, it is the first time that the slip with friction conditions are imposed in this context. 
This extends the Navier boundary condition to the interface case and interpolates between no-slip of the bulk phase and embedded material layer, as considered in
\cite{barrett2014stable,BothePruess2010,nitschke_voigt_wensch_2012,ReuskenZhang}, and  uncoupled lateral dynamics 
of fluidic surfaces as studied, for example, in \cite{reuther2015interplay,rodrigues2015semi,torres2019modelling}. The model contributes to  characterizing the interaction between a rather viscous and dense interface material and  less viscous and dense bulk fluid,
as is the case of, e.g., a lipid vesicle or cell
membrane with its content and surrounding extracellular fluid. 

We first state the problem in terms of incompressible two-phase Navier--Stokes flow in the bulk coupled to
surface fluid equations and show the energy balance for this system. This balance
demonstrates that the system with the proposed coupling conditions is thermodynamically consistent, i.e.~in the absence of external forces 
and with no energy inflow through the boundary the system is dissipative.
Then, we introduce a simplifying assumption: the coupled system 
has reached a steady state and inertia terms can be neglected. Although the resulting surface--bulk Stokes problem
is a strong simplification of the original problem, it still retains interesting features for the purpose of numerical analysis.

For the simplified model, we show well-posedness following the abstract  framework for
generalized saddle point problems proposed in, e.g., \cite{bernardi1988generalized,nicolaides1982existence}.
Next, we present a sharp interface, geometrically unfitted finite element (FE) method.
Unfitted methods allow the sharp interface to cut through the elements of a  fixed background grid. 
Their main advantage is the relative ease in handling time-dependent domains, implicitly defined interfaces, 
and problems with strong geometric deformations \cite{Bordas2018}.
For the bulk two-phase Stokes problem, we apply an isoparametric unfitted finite element approach~\cite{lehrenfeld2016high} of the CutFEM (or Nitsche-XFEM) family\cite{CLAUS2019185,FRACHON201977,HANSBO201490,He_Song2019,Massing_Larson2014,Wang_Chen2019}. For a (bulk only) interface Stokes problem with slip between phases, the approach was studied in \cite{olshanskii2021unfitted}.
CutFEM uses overlapping fictitious domains in combination with ghost
penalty stabilization \cite{B10} to enrich and stabilize the solution.
In this paper, we consider the unfitted generalized Taylor--Hood finite element pair
$\bP_{k+1}-P_k$, $k\ge1$. For more details on the isoparametric unfitted finite element, we refer to
\cite{Lehrenfeld2017,lehrenfeld2018analysis}.
For the discretization of the surface Stokes problem, we apply the trace finite element method \cite{Olshanskii_Reusken2019,jankuhn2020error},
which uses traces of the bulk finite element functions. In a set of numerical experiments, we address
the dependence of the discretization error of the method on variations in the friction coefficient, ratio of bulk fluid viscosity, surface fluid viscosity, mesh size, and position of the interface relative to the fixed computational mesh.

One option to solve the coupled bulk-surface problem is to adopt 
a monolithic approach. However, such an approach would quickly lead to
high computational costs as the mesh gets refined. Instead, we introduce a partitioned scheme based on fixed-point iterations  \cite{B-quarteroniv1}. With this algorithm, 
the bulk and surface flow problems are solved separately and sequentially
and the coupling conditions are enforced in an iterative fashion. We choose this algorithm for its
simplicity of implementation and study its performance when model and discretization parameters vary. 
Our numerical experiments provide evidence of the robustness of the proposed approach with respect to the contrast in viscosity in the bulk fluid, surface fluid viscosity, and position of the interface relative to the background mesh. At the same time,  convergence slowdown was observed for certain values of the slip coefficients. 

The outline of the paper is as follows.
In Sec.~\ref{sec:t_dep}, we introduce the strong formulation of the coupled problem
and the associated energy balance. Sec.~\ref{sec:simplified_p} presents 
the strong and weak formulations of the simplified problem, together with the finite element discretization.
In Sec.~\ref{sec:part_meth}, we propose a partitioned algorithm for the numerical
solution of the coupled problem. 
Numerical results in 3 dimensions are reported in Sec.~\ref{sec:num_res}.
For all the simulations in this paper, we have used NGsolve \cite{NGSolve,Gangl2020},
a high performance multiphysics finite element software with a  Python interface,
and add-on library ngsxfem \cite{ngsxfem}.

\section{A two-phase fluid with material viscous interface}\label{sec:t_dep}
\begin{wrapfigure}{r}{0.49\textwidth}\label{Geom}
    \centering
    \includegraphics[scale=0.25]{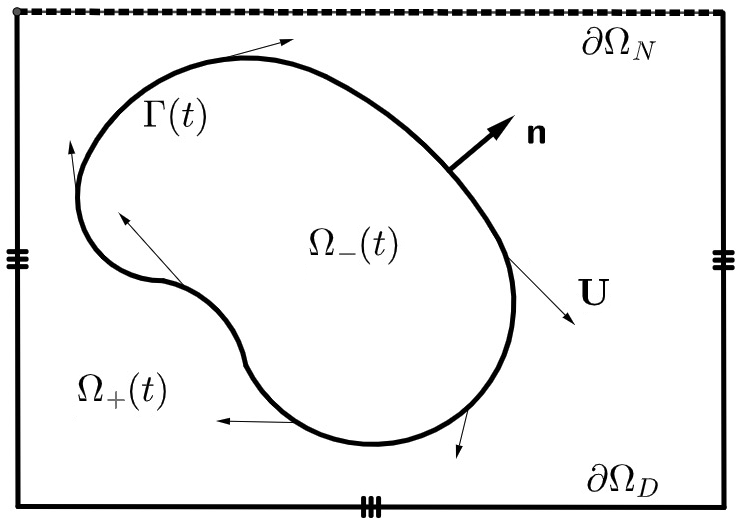}
    \caption{Fluid domains and interface in $\mathbb{R}^2$. 
    }
    \vskip-3ex
    \label{fig:Geom}
\end{wrapfigure}

Consider a fixed volume $\Omega\subset \mathbb{R}^3$ filled with two immiscible, viscous,
and incompressible fluids separated by  an interface  $\Gamma(t)$ for all $t \in [0,T]$.
We assume that $\Gamma(t)$ stays closed and sufficiently smooth (at least $C^2$) for all $t \in [0,T]$. Surface $\Gamma(t)$ separates $\Omega$ into
two phases  (subdomains) $\Omega_+(t)$ and $\Omega_-(t):=\Omega\setminus \overline{\Omega_{+}(t)}$.
We assume $\Omega_-(t)$ to be completely internal, i.e.~$\partial\Omega_-(t)\cap \partial\Omega = \emptyset$ for all times.
See Fig.~\ref{fig:Geom}.
Denote by $\bn^{\pm}$ the outward normals for $\Omega_{\pm}(t)$
and $\bn$ the normal on $\Gamma$ pointing from $\Omega_-(t)$ to $\Omega_+(t)$: it holds that $\bn^- = \bn$ and $\bn^+ = - \bn$ at $\Gamma$.
For ease of notation, from now on we will drop the dependance on $t$ for $\Gamma$, $\Omega_+$,
and $\Omega_-$.

To describe the fluid motion, consider the bulk fluid velocities $\bu^{\pm}(t):\Omega_{\pm}\to \mathbb{R}^3$ and pressures $p^{\pm}(t):\Omega_{\pm}\to \mathbb{R}$.
The motion of the fluids occupying subdomains $\Omega_{\pm}$ is governed by the incompressible Navier--Stokes equations
\begin{align}
    \rho^{\pm}\dot{\partial}_t \bu^{\pm}&=\Div \bsigma^{\pm} + \bbf^{\pm} &\text{ in }\Omega_{\pm}, \label{eq:NS1} \\
    \Div {}\bu^{\pm}&=0 &\text{ in }\Omega_{\pm}, \label{eq:NS2}
\end{align}
for all $t\in(0,T)$.
In \eqref{eq:NS1}, constants $\rho^{\pm}$ represent the fluid density, $\dot{\partial}_t$ denotes the material derivative,
$\bbf^{\pm}$ are the external body forces, and $\bsigma^{\pm}$ are the Cauchy stress tensors.
For Newtonian fluids in both phases, the Cauchy stress tensors have the following expression:
\begin{align}\label{eq:sigmas}
\bsigma^{\pm}= -p^{\pm} \mathbf{I} +2\mu^{\pm} \bD(\bu^{\pm}), \quad \bD(\bu^{\pm}) = \frac{1}{2} (\nabla \bu^{\pm} + (\nabla\bu^{\pm})^T) \text{ in }\Omega_{\pm},
\end{align}
where constants $\mu^{\pm}$ represent the fluid dynamic viscosity.

We assume interface $\Gamma$ to be a thin  material layer with possibly different material properties
from the bulk fluid. Motivated by applications in cell biology, we consider a viscous inextensible interface modeled
as an ``incompressible'' surface  fluid.
 The evolution of the material interface can be described in terms of the velocity of this surface fluid denoted by $\bU$.
 Later, we will need the decomposition of $\bU$ into tangential and normal components:
 $\bU=\bU_T+U_N\bn$, with $\bU_T\cdot\bn=0$, $U_N=\bn\cdot\bU$.
The surface Navier-Stokes equations governing the motion of a fluidic deformable layer appear in several work~ \cite{yavari2016nonlinear,salbreux2017mechanics,jankuhn2017incompressible,miura2018singular}. Here, we adopt the formulation in terms of tangential differential operators from \cite{jankuhn2017incompressible}, where the equations  have been derived
from conservation principles.
To introduce these equations, we need some further definitions.
Let  $\bP(\bx):=\bI-\bn(\bx)\bn(\bx)^T$ for $\bx \in \Gamma$ be the orthogonal projection onto the tangent plane.
  For a scalar function $\pi:\, \Gamma \to \mathbb{R}$ or a vector function $\bU:\, \Gamma \to \mathbb{R}^3$  we define $\pi^e\,:\,\mathcal{O}(\Gamma)\to\mathbb{R}$, $\bU^e\,:\,\mathcal{O}(\Gamma)\to\mathbb{R}^3$, smooth  extensions of $\pi$ and $\bU$ from $\Gamma$ to its neighborhood $\mathcal{O}(\Gamma)$.
 The surface gradient and covariant derivatives on $\Gamma$ are then defined as $\nablaG \pi=\bP\nabla \pi^e$ and  $\nabla_\Gamma \bU:= \bP \nabla \bU^e \bP$. The definitions  of $\nablaG \pi$ and $\nabla_\Gamma \bU$ are  independent of the particular smooth extension of $\pi$ and $\bU$ off $\Gamma$.
On $\Gamma$, we consider the surface rate-of-strain tensor \cite{GurtinMurdoch75} given by
\begin{equation} \label{strain}
 \bD_\Gamma(\bU):= \frac12 \bP (\nabla \bU +(\nabla \bU)^T)\bP = \frac12(\nabla_\Gamma \bU + (\nabla_\Gamma \bU)^T).
 \end{equation}
The surface divergence for a vector $\bg: \Gamma \to \R^3$ and
a tensor $\bA: \Gamma \to \mathbb{R}^{3\times 3}$ are defined as:
$
 \divG \bg  := \tr (\gradG \bg)$, $
 \divG \bA  := \mbox{\small $\left( \divG (\be_1^T \bA),\,
               \divG (\be_2^T \bA),\,
               \divG (\be_3^T \bA)\right)^T$ },
 $
 where $\be_i$ is the $i$th standard basis vector.
With this notation, the surface Navier--Stokes equations take the form:
   \begin{align}
 \rho_{\Gamma}\dot{\partial}_t\bU&=-\nabla_{\Gamma}\pi+2\mu_{\Gamma}\divG \bD_\Gamma(\bU)+\bbf_\Gamma+\bb^e+\pi\kappa\bn
 &\text{ on }\Gamma, \label{eq:sNS1}\\
         \divG \bU &=0 &\text{ on }\Gamma, \label{eq:sNS2}
   \end{align}
 where $ \rho_{\Gamma}$ is the surface fluid density, $\mu_{\Gamma}$ is the surface fluid dynamic viscosity,
$\kappa$ denotes  point-wise doubled mean curvature on $\Gamma$, and $\pi$ is the surface fluid pressure.
 The material derivative in \eqref{eq:sNS1} is taken with respect to surface fluid trajectories, i.e. $\dot{\partial}_t\bU= \frac{\partial \bU}{\partial t}+(\bU\cdot\nabla)\bU$. Note that $\dot{\partial}_t\bU$ is an intrinsic surface quantity, although both terms $\frac{\partial \bU}{\partial t}$ and $(\bU\cdot\nabla)\bU$ depend on extension of $\bU$ in the bulk.  On the right hand side of \eqref{eq:sNS1},
 $\bbf_\Gamma$ denotes the external area force acting on the  surface as a result of the interaction with the bulk fluids
(specified below),  while $\bb^e$ denotes other possible area force
(such as elastic bending forces) and not specified further in the present paper.

Next, we turn to the coupling conditions between equations \eqref{eq:NS1}--\eqref{eq:NS2} posed in the bulk
 and equations \eqref{eq:sNS1}--\eqref{eq:sNS2} posed on $\Gamma$.
 First, the immiscibility condition means that the bulk fluid does not penetrate through $\Gamma$, which  implies that
 \begin{equation}
     \bu^+\cdot \bn= U_N=\bu^-\cdot \bn\qquad \qquad \text{ on }\Gamma. \label{eq:cc1}
 \end{equation}
Normal velocity $U_N$ determines radial deformations of $\Gamma(t)$ and so it governs the geometric evolution of the interface, which
can be defined through the Lagrangian mapping
$\Psi(t, \cdot)$ from $\Gamma(0)$ to $\Gamma(t)$: for $\bx\in\Gamma(0)$, $\Psi(t,\bx)$ solves the ODE system
\begin{equation}\label{Lagrange}
\Psi(0,\bx)=\bx,\quad  \frac{\partial \Psi(t,\bx)}{\partial t}=U_N(t,\Psi(t,\bx)),\quad t\in[0,T].
\end{equation}
In fluid vesicles and cells, typically a viscous and dense lipid membrane, represented by $\Gamma$, is
surrounded by a less viscous and less dense liquid. We are interested in modelling
slip with friction between the bulk fluid and the viscous  membrane.
Thus, we consider Navier-type conditions
   \begin{align}
  \bP{\bsigma^+\bn}&=f^+(\bP\bu^+-\bU_T)  &\text{ on }\Gamma,  \label{eq:cc2}\\
  \bP{\bsigma^-\bn}&=-f^-(\bP\bu^--\bU_T)  &\text{ on }\Gamma,  \label{eq:cc3}
     \end{align}
  where $f^-$ and $f^+$ are friction coefficients at $\Gamma$ on the $\Omega_-$ and $\Omega_+$ side, respectively.
  Conditions~\eqref{eq:cc2}--\eqref{eq:cc3} model an incomplete adhesion of a bulk fluid to the material surface with $1/f^\pm$ often referred to as a ``slip length''~\cite{navier1823memoire}; see, e.g., \cite{bocquet2007flow,lauga2007microfluidics} for the modern description of experimental and theoretical validations.  In particular, the acceptance of non-zero slip length resolves the well-known ``no-collision paradox'' ~\cite{cooley1969slow,hocking1973effect,gerard2015influence}, thus suggesting \eqref{eq:cc2}--\eqref{eq:cc3} to be an important modeling assumption in the simulation of a lipid vesicle -- cell membrane contact (and fusion). We finally note that the Navier conditions are not an uncommon choice in numerical models, if the flow in  boundary region is under-resolved~\cite{john2002slip}. 
  
The area force in \eqref{eq:sNS1} coming from the bulk fluid is defined by the jump of the normal stress on $\Gamma$:
\begin{equation}
  \bbf_\Gamma = [\bsigma\bn]^+_-= \bsigma^+\bn - \bsigma^-\bn ~ \text{ on }\Gamma.  \label{eq:cc4}
\end{equation}

We summarize the complete system of equations and coupling conditions below
\begin{equation}\label{eq:Full}
\left\{
\begin{aligned}
      \rho^{\pm}\dot{\partial}_t \bu^{\pm}-\mu^{\pm}\Delta\bu^\pm+\nabla p^{\pm}&=  \bbf^{\pm} &\text{ in }\Omega_{\pm}, \\
       \Div {}\bu^{\pm}&=0 &\text{ in }\Omega_{\pm},\\
       \rho_{\Gamma}\dot{\partial}_t\bU-2\mu_{\Gamma}\divG \bD_\Gamma(\bU)+\nabla_{\Gamma}\pi-\pi\kappa\bn&=[\bsigma\bn]^+_- +\bb^e
 &\text{ on }\Gamma, \\
         \divG \bU &=0 &\text{ on }\Gamma, \\
         \bu^+\cdot \bn=\bu^-\cdot \bn& = U_N& \text{ on }\Gamma,\\
         \bP{\bsigma^\pm\bn}&=\pm f^\pm(\bP\bu^\pm-\bU_T)  &\text{ on }\Gamma,
    \end{aligned}
    \right.
\end{equation}
with the evolution of  $\Omega_{\pm}$ and $\Gamma$ defined by the velocity solving the system through \eqref{Lagrange}.
On $\partial\Omega$ the system is endowed with boundary conditions either for the bulk velocity or for the bulk  normal stress:
\begin{align}
    \bu^+&=\bg\quad \text{ on }\partial\Omega_D, \label{eq:bcD} \\
    \bsigma^+\bn^+ &= \bbf_N\quad \text{ on }\partial\Omega_N. \label{eq:bcN}
\end{align}
 Here $\overline{\partial\Omega_D}\cup\overline{\partial\Omega_N}=\overline{\partial\Omega}$ and $\partial\Omega_D \cap\partial\Omega_N=\emptyset$. See Fig.~\ref{fig:Geom}. At $t = 0$, initial velocity is given $\bu^{\pm} = \bu_0^{\pm}$ in $\Omega_{\pm}(0)$
 and $\bU = \bU_0$ on $\Gamma(0)$.

\subsection{Energy balance of the continuous coupled problem}
We look for the energy balance of the coupled system \eqref{eq:Full}--\eqref{eq:bcN}.
We make use of the following identities for time-dependent domains $\Omega_\pm(t)$, whose
only moving part of the boundary is $\Gamma(t)$:
\begin{align}
\frac{1}{2}\frac{d}{dt}\int_{\Omega_{\pm}(t)} |\bu^{\pm}|^2\,dV &= \int_{\Omega_{\pm}(t)} \bu^{\pm}\cdot{}\frac{\partial \bu^{\pm}}{\partial t}\,dV +
\frac{1}{2}\int_{\Gamma(t)} |\bu^{\pm}|^2 \bu^{\pm} \cdot \bn^{\pm} \,dS, \label{eq:trasport}\\
\int_{\Omega_{\pm}(t)} (\bu^{\pm} \cdot \nabla) \bu^{\pm} \cdot \bu^{\pm} \,dV &= \frac{1}{2}\int_{\partial\Omega_{\pm}(t)} |\bu^{\pm}|^2 \bu^{\pm} \cdot \bn^{\pm} \,dS. \label{eq:conv_i}
\end{align}
Identity \eqref{eq:trasport} is the Reynolds transport theorem, while identity \eqref{eq:conv_i} is obtained from integration by parts.

 Let us start from the kinetic energy of the fluid in $\Omega_-$:
\begin{align}
        \frac{d}{dt}E^-&=\frac{1}{2}\frac{d}{dt}\int_{\Omega_-}\!\!\rho^- |\bu^-|^2\,dV =  \int_{\Omega_-}\rho^- \bu^-\cdot{}\frac{\partial \bu^-}{\partial t}\,dV + \frac{1}{2}\int_{\Gamma} \rho^- |\bu^-|^2 \bu^- \cdot \bn \,dS \cl
        &=\int_{\Omega_-} \bu^-\cdot(\Div \bsigma^{-} + \bbf^- - \rho^- (\bu^- \cdot \nabla) \bu^-) \,dV + \frac{1}{2}\int_{\Gamma} \rho^- |\bu^-|^2 \bu^- \cdot \bn \,dS\cl
       &= \int_{\Omega_-}\nabla \cdot(\bsigma^- \bu^-)\,dV -\int_{\Omega_-}\bsigma^-:\nabla\bu^-\,dV + \int_{\Omega_-} \bu^-\cdot \bbf^-\,dV \cl
 &=\int_\Gamma \bu^-\cdot (\bsigma{}^-\bn)\,dS -\int_{\Omega_-} \bsigma^- : \bD(\bu^-)dV + \int_{\Omega_-} \bu^-\cdot \bbf^-\,dV\cl
 &=\int_\Gamma \bu^-\cdot (\bsigma{}^-\bn)\,dS -2\mu^{-}\int_{\Omega_-} \| \bD(\bu^-) \|^2\,dV + \int_{\Omega_-} \bu^-\cdot \bbf^-\,dV \label{eq:E-}
\end{align}
Above, we have used \eqref{eq:trasport}, \eqref{eq:NS1}, \eqref{eq:NS2}, \eqref{eq:conv_i}, and integration by parts.

We repeat similar steps for the kinetic energy of the fluid in $\Omega_+$, the main difference being that
$\overline{\partial\Omega_+} = \overline{\Gamma} \cup \overline{\partial\Omega_D}\cup\overline{\partial\Omega_N}$
while $\overline{\partial\Omega_-}= \overline{\Gamma}$.
We obtain:
\begin{align}
        &\frac{d}{dt}E^+=- \int_{\Gamma}\bu^+\cdot{}(\bsigma^+\bn)\,dS -2\mu^{+}\int_{\Omega_+}\norm{\bD(\bu^+)}^2 dV + \int_{\Omega_+} \bu^+\cdot \bbf^+\,dV +B \label{eq:E+}
\end{align}
where
\begin{align}
B=\int_{\pa\Omega_D}\left( \bg\cdot (\bsigma^+\bn) - \frac{1}{2} \rho^+ |\bg|^2 \bg \cdot \bn^+ \right)\,dS+
\int_{\pa\Omega_N} \left( \bu^+\cdot{}\bbf{} - \frac{1}{2} \rho^+ |\bu^+|^2 \bu^+ \cdot \bn^+ \right)\,dS. \label{eq:B}
\end{align}

Putting together \eqref{eq:E-} and \eqref{eq:E+}, we obtain the total kinetic energy for the bulk flow:
 \begin{align}
        \frac{dE}{dt}=&\frac{d}{dt}(E^++E^-)=-2\mu^-\int_{\Omega_-} \| \bD(\bu^-) \|^2\,dV
        -2\mu^+\int_{\Omega_+}\norm{\bD(\bu^+)}^2 dV -\int_\Gamma U_N[\bn^T\sigma\bn]^{+}_{-}dS \cl
        &-\int_{\Gamma}f^-\bP\bu^-\cdot(\bP\bu^--\bU_T)dS-\int_{\Gamma}f^+\bP\bu^+\cdot(\bP\bu^+-\bU_T)dS \cl
        &+ \int_{\Omega^{-}} \bu^-\cdot \bbf^-\,dV + \int_{\Omega^{+}} \bu^+\cdot \bbf^+\,dV +B \label{eq:E}
\end{align}
where we have used \eqref{eq:cc1}, \eqref{eq:cc2}, and \eqref{eq:cc3}.
%
%
For the kinetic energy of the surface flow, we will use
the surface analogue of the Reynolds transport theorem (see, e.g. \cite[Lemma 2.1]{DziukElliot2013a}):
\begin{align}\label{eq:transport_rhoxi}
\frac{{\rm d}}{{\rm d}t} \int_{\Gamma} f  \, dS =  \int_{\Gamma} \left(\dot{\partial}_t f + f \divG \bU \right) \, dS.
\end{align}
We also need the surface integration by part identity:
\[
\int_\Gamma(g\divG\blf+ \blf\cdot\nablaG g)\,dS= \int_\Gamma \kappa g(\blf\cdot\bn)\,dS,
\]
which is valid for a smooth closed surface with smooth scalar field $g$ and vector field $\blf$.
A componentwise application of this equality yields the identity
\begin{equation}
    \int_{\Gamma} \bg\cdot\divG\bG\,dS=-\int_{\Gamma} \bG:\nabla_{\Gamma}\bg\, dS\label{eq:lemma1}
\end{equation}
for a vector field $\bg\in$ $({C}^1(\Gamma))^3$ and a matrix function $\bG\in ({C}^1(\Gamma))^{3\times{}3}$ such that $\bG=\bP\bG\bP$.

The energy balance on $\Gamma$ is given by:
\begin{align}
    \frac{dE_\Gamma}{dt}&= \frac{1}{2}\frac{d}{dt}\int_{\Gamma}\!\!\rho_\Gamma |\bU|^2\,dS=\int_{\Gamma}\left(\rho_\Gamma\bU\cdot \dot{\partial}_t\bU + \frac{1}{2} \rho_\Gamma |\bU|^2 \divG \bU \right) \,dS
    \cl
        &=\int_{\Gamma} \bU\cdot \left(-\nabla_{\Gamma}\pi+2\mu_{\Gamma}\divG \bD_\Gamma(\bU)+\bbf_\Gamma+\bb^e\right)\,dS  \cl
        &=   - 2 \mu_\Gamma \int_{\Gamma} \bD_\Gamma(\bU) : \nabla_\Gamma \bU  \, dS  +
        \int_{\Gamma} \bU\cdot  [\bsigma\bn]^+_- \,dS + \int_{\Gamma} \bU_T\cdot \bb^e\,dS  \cl
                &=  - 2 \mu_\Gamma \int_{\Gamma} \| \bD_\Gamma(\bU) \|^2  \, dS +  \int_{\Gamma} U_N  [\bn^T\bsigma\bn]^+_- \,dS  \cl
        &\quad+ \int_{\Gamma} \bU_T\cdot  \left(f^+(\bP\bu^+-\bU_T) +f^-(\bP\bu^--\bU_T) \right) \,dS + \int_{\Gamma} \bU\cdot \bb^e \,dS,  \label{eq:E_gamma}
     \end{align}
where we have applied \eqref{eq:lemma1} and used \eqref{eq:transport_rhoxi}, \eqref{eq:sNS1}-\eqref{eq:cc1},
\eqref{eq:cc2}-\eqref{eq:cc4}, and integration by parts.

Combining \eqref{eq:E} and \eqref{eq:E_gamma}, we get the kinetic energy for the bulk and surface flows:
\begin{align}
    \frac{dE}{dt}+\frac{dE_{\Gamma}}{dt}&=
      \underbrace{-2\mu^{-}\int_{\Omega_-} \| \bD(\bu^-) \|^2\,dV -2\mu^+\int_{\Omega_+}\norm{\bD(\bu^+)}^2 dV}_{\text{Bulk fluid viscous dissipation}}
        \underbrace{-2 \mu_\Gamma \int_{\Gamma} \| \bD_\Gamma(\bU_T) \|^2  \, dS}_{\text{Surface fluid viscous dissipation}}
        \cl
        & \quad \underbrace{-\int_{\Gamma}f^-\|\bP\bu^--\bU_T\|^2dS-\int_{\Gamma}f^+\|\bP\bu^+-\bU_T\|^2d} _{\text{Frictional energy dissipation}} \cl
        & \quad \underbrace{+ \int_\Gamma \bu^-\cdot \bbf^-\,dV + \int_\Gamma \bu^+\cdot \bbf^+\,dV+ \int_{\Gamma} \bU\cdot \bb^e dS}_{\text{work of external forces }}   \underbrace{\phantom{\int_\Gamma}\hspace{-3ex}+B}_{\text{work of b.c.}} \label{eq:E_tot}
\end{align}
where $B$, i.e.~the work of the boundary conditions, is defined in \eqref{eq:B}.

From \eqref{eq:E_tot}, we see that in the absence of external forces and with no  energy inflow through the 
boundary the system is dissipative, i.e. thermodynamically consistent.

\section{A simplified steady problem}\label{sec:simplified_p}
In this section, we consider a (strongly) simplified version of the problem presented in Sec.~\ref{sec:t_dep}. Our main assumption is that the coupled bulk and surface fluid system has reached a steady state and inertia terms can be neglected. Since the steady state implies $\Gamma(t)=\Gamma(0)$, we have $U_N=0$ and hence $\bU=\bU_T$.
This simplified surface--bulk Stokes problem models a viscosity dominated two-phase flow with the viscous interface in a dynamical equilibrium; see also Remark~\ref{Rem}. It is an interesting model problem for the purpose of numerical analysis.
With these simplifications, the equations \eqref{eq:NS1}--\eqref{eq:NS2} become:
\begin{align}
    -\mu^\pm\Delta\bu^{\pm}+\nabla p &= \bbf^{\pm} \qquad\text{in }\Omega_{\pm}, \label{eq:Stokes1} \\
    \Div {}\bu^{\pm}&=0 \qquad ~~ \text{in }\Omega_{\pm}. \label{eq:Stokes2}
\end{align}
We impose a non-homogeneous Dirichlet condition on the entire outer boundary of $\Omega$, i.e.~problem
\eqref{eq:Stokes1}--\eqref{eq:Stokes2} is supplemented with boundary condition \eqref{eq:bcD} with $\bg\in [H^{1/2}(\pa\Omega_D)]^3$ and
$\partial\Omega_D = \partial\Omega$.
Under our assumption, the momentum equation for the surface fluid simplifies to
$
-2\mu_{\Gamma}\divG \bD_\Gamma(\bU_T)+\nabla_{\Gamma}\pi-\pi\kappa\bn=[\bsigma\bn]^+_- +\bb^e.
$
The tangential part of the above momentum equation together with the inextensibility condition \eqref{eq:sNS2} leads to the surface Stokes problem
   \begin{align}
-2\mu_{\Gamma}\bP\divG \bD_\Gamma(\bU_T)+\nabla_{\Gamma}\pi&= [\bP\bsigma\bn]^+_-+ \bP\bb^e\qquad
 &\text{ on }\Gamma, \label{eq:sStokes1}\\
         \divG \bU_T &=0\qquad&\text{ on }\Gamma, \label{eq:sStokes2}
   \end{align}
while the normal part  simplifies to
\begin{equation}
  [\bn^T\bsigma\bn]^-_+=\pi\kappa\qquad\text{ on }\Gamma. \label{eq:scc4b}
\end{equation}
The  interface condition above is standard in many models of two-phase flows, where $\pi$ has the meaning of the surface tension coefficient.

Coupling condition \eqref{eq:cc1} is replaced by:
\begin{equation}
    \bu^+\cdot \bn= \bu^-\cdot \bn\qquad \text{ on }\Gamma, \label{eq:scc1}
\end{equation}
while  conditions \eqref{eq:cc2} and \eqref{eq:cc3} still hold:
\begin{equation}\label{eq:scc4a}
\bP{\bsigma^\pm\bn}=\pm f^\pm(\bP\bu^\pm-\bU_T)\qquad\text{ on }\Gamma.
\end{equation}
Finally, we will see that the (weak formulation of the) problem is well-posed under two mean conditions for the bulk pressure:
\[
\int_{\Omega^\pm} p^\pm\,dx=0.
\]

\begin{remark}\label{Rem}\rm Since $\bU\cdot\bn=0$,  condition \eqref{eq:scc1} allows the flow
through the steady interface $\Gamma$. This is inconsistent with \eqref{eq:cc1}, which assumes
immiscibility of fluids. For a physically consistent formulation that describes the true  equilibrium one has to set $\bu^+\cdot \bn= \bu^-\cdot \bn=0$ on  $\Gamma$, but allow the shape of  $\Gamma$ to be the unknown, i.e.~to be determined as a part of the problem. For such equilibrium to exist, external forces and boundary conditions may have to satisfy additional constraints. Finding such constraints
and solving the resulting non-linear problem is outside the scope of this paper. We rather follow a common convention in the analysis
of models for steady two-phase problems and allow \eqref{eq:scc1} for the steady interface; see, e.g.~\cite{Ganesan05,Olshanskii04,gross2011numerical,bootland2019preconditioners}.
\end{remark}

\subsection{Variational formulation}

The purpose of this section is to derive the variational formulation of coupled problem \eqref{eq:Stokes1}--\eqref{eq:scc4a}.
Let us introduce some standard notation. The space of functions whose square is integrable in a domain $\omega$
is denoted by $L^2(\omega)$.
The space of functions whose distributional derivatives of order up to $m \geq 0$  (integer)
belong to $L^2(\omega)$ is denoted by $H^m(\omega)$.
The space of vector-valued functions with components in $L^2(\omega)$ is denoted with $L^2(\omega)^3$.
$H^1(\div,\omega)$ is the space of functions in $L^2(\omega)$ with divergence in $L^2(\omega)$.
Moreover, we introduce the following functional spaces:
\begin{align}
V^- &= H^1(\Omega_-)^3, ~V^+ =\{\bu \in H^1(\Omega_+)^3, \bu{\big|_{\partial \Omega_D}}=\bg \}, ~V^+_0 =\{\bu \in H^1(\Omega_+)^3, \bu{\big|_{\partial \Omega_D}}=\boldsymbol{0} \}, \cl
V^\pm & = \{ \bu = (\bu^-,\bu^+) \in V^- \times V^+ , \bu^- \cdot \bn = \bu^+ \cdot \bn~\text{on}~\Gamma \}, \cl
V^\pm_0 & = \{ \bu = (\bu^-,\bu^+) \in V^- \times V^+_0 , \bu^- \cdot \bn = \bu^+ \cdot \bn~\text{on}~\Gamma \}, \cl
L^2_{\pm} & = \{ p = (p^-,p^+) \in L^2(\Omega_-)\times L^2(\Omega_+),~s.t.~\int_{\Omega^\pm} p^\pm\,dx=0 \}, \cl
V_\Gamma&=\{\bU\in H^1(\Gamma)^3 :\bU\cdot\bn=0\}. \el
\end{align}
The space $V^\pm$ can be also characterized as $(V^- \times V^+)\cap H^1(\div,\Omega)$.
We use $( \cdot, \cdot)_\omega$ and $\langle\cdot, \cdot\rangle_\omega$ to denote the $L^2$ product and the duality pairing, respectively.

Multiplying  \eqref{eq:Stokes1} by $\bv\in V^{\pm}_0$ and \eqref{eq:Stokes2} by $q \in L^2_0(\Omega)$ and integrating over each subdomain, we see that smooth bulk velocity and pressure satisfy integral identity:
\begin{align}
&- (p^-,\Div  \bv^-)_{\Omega_{-}} - (p^+,\Div  \bv^+)_{\Omega_{+}} + 2(\mu^-\bD(\bu^-),\bD(\bv^-))_{\Omega_{-}}  + 2(\mu^+\bD(\bu^+),\bD(\bv^+))_{\Omega_{+}}  \cl
&\quad \quad  - \langle \pi \kappa , \bv^- \cdot \bn \rangle_\Gamma + \langle f^-(\bP\bu^--\bU ), \bP\bv^- \rangle_\Gamma + \langle f^+(\bP\bu^+-\bU ), \bP\bv^+ \rangle_\Gamma   \cl
&\quad \quad = (\bbf^-,\bv^-)_{\Omega_{-}} + (\bbf^+,\bv^+)_{\Omega_{+}}  \label{eq:waek-1} \\
&(\Div \bu^{-}, q^-)_{\Omega_{-}}+ (\Div \bu^{+}, q^+)_{\Omega_{+}}=0 \label{eq:waek-2}
\end{align}
for all $(\bv, q) \in V^{\pm}_0 \times L^2_0(\Omega)$.
The interface terms in \eqref{eq:waek-1} have been obtained using
coupling conditions \eqref{eq:scc4a} and \eqref{eq:scc4b} as follows:
\begin{align}
-\langle \bsigma^- \bn, \bv^- \rangle_\Gamma  + \langle \bsigma^+ \bn, \bv^+\rangle_\Gamma &=
-\langle \bP\bsigma^- \bn, \bP\bv^- \rangle_\Gamma + \langle \bP\bsigma^+ \bn, \bP\bv^+\rangle_\Gamma - \langle [\bn^T\bsigma \bn]^-_+, \bv^- \cdot \bn \rangle_\Gamma \cl
& = \langle f^-(\bP\bu^--\bU ), \bP\bv^- \rangle_\Gamma + \langle f^+(\bP\bu^+-\bU ), \bP\bv^+ \rangle_\Gamma \cl
& \quad  - \langle \pi \kappa , \bv^- \cdot \bn \rangle_\Gamma. \el
\end{align}
Likewise, we find that the surface velocity and pressure satisfy the following integral identities:
\begin{align}
  &-(\pi,\divG  \bV)_{\Gamma}+2(\mu_{\Gamma}\bD_{\Gamma}(\bU ),\bD_{\Gamma}(\bV))_{\Gamma}
  - \langle f^-(\bP\bu^--\bU ), \bV \rangle_\Gamma \cl
  & \quad \quad- \langle f^+(\bP\bu^+-\bU ), \bV \rangle_\Gamma=(\bP\bb^e,\bV)_{\Gamma}\label{eq:sStokes_weak1}\\
&  (\divG \bU , \tau)_\Gamma =0  \label{eq:sStokes_weak2}
\end{align}
for all $(\bV, \tau) \in V_\Gamma \times L^2_0(\Gamma)$. 

The weak formulation of the coupled problem \eqref{eq:Stokes1}--\eqref{eq:scc4a} follows by combining \eqref{eq:waek-1}--\eqref{eq:waek-2} and \eqref{eq:sStokes_weak1}--\eqref{eq:sStokes_weak2}. In order to write it,
we introduce the following forms for all $\bu \in V^{\pm}$, $\bv\in\bV^\pm_0$, $\bU, \bV\in V_\Gamma$, $p \in L^2(\Omega)$, $\pi\in L^2(\Gamma)$:
\begin{align}
    a(\{\bu,\bU\},\{\bv,\bV\})=& 2(\mu^-\bD(\bu^-),\bD(\bv^-))_{\Omega_{-}}  + 2(\mu^+\bD(\bu^+),\bD(\bv^+))_{\Omega_{+}} \cl
    & +2(\mu_{\Gamma}\bD_{\Gamma}(\bU ),\bD_{\Gamma}(\bV))_{\Gamma} + \langle f^-(\bP\bu^--\bU ), \bP\bv^- - \bV \rangle_\Gamma \cl
    & + \langle f^+(\bP\bu^+-\bU ), \bP\bv^+ - \bV  \rangle_\Gamma, \cl
    b(\{\bv,\bV\},\{p,\pi\})=&- (p^-,\Div  \bv^-)_{\Omega_{-}} - (p^+,\Div  \bv^+)_{\Omega_{+}}-(\pi,\divG  \bV)_{\Gamma},\cl
    s(\bv,\pi)=&-\langle \pi \kappa , \bv \cdot \bn \rangle_\Gamma,\cl
    r(\bv,\bV)=&(\bbf^-,\bv^-)_{\Omega_{-}} + (\bbf^+,\bv^+)_{\Omega_{+}}+(\bP\bb^e,\bV)_{\Gamma}. \el
\end{align}
Then, the weak formulation reads:
{\it Find $(\bu, p) \in V^\pm \times L^2_\pm$,
and $(\bU , \pi) \in V_\Gamma \times L^2(\Gamma)$ such that
\begin{equation}\label{weak}
\left\{
\begin{split}
   a(\{\bu,\bU \},\{\bv,\bV\})+b(\{\bv,\bV\},\{p,\pi\})+s(\bv,\pi)&= r(\bv,\bV) \\
   b(\{\bu,\bU \},\{q,\tau\})&=0
\end{split}
\right.
\end{equation}
for all $(\bv, q) \in V_0^{\pm} \times L^2_0(\Omega)$ and $(\bV, \tau) \in V_\Gamma \times L^2_0(\Gamma)$.} Note that test and trial pressure spaces both involve two (different) gauge conditions.

\subsection{Well-Posedness}\label{sec:well}
With the goal of proving the well-posedness of the stationary problem,
we start by showing that $a(\{\cdot,\cdot\},\{\cdot,\cdot\})$ is coercive.
Let $\norm{\bu}^2_{H^1(\Omega_{\pm})}=\norm{\bu^+}^2_{H^1(\Omega_+)}+\norm{\bu^-}^2_{H^1(\Omega_-)}$.
Let $\norm{p}_{L^2(\Omega_{\pm})}^2=\norm{p }_{\Omega_{+}}^2+\norm{p }_{\Omega_{-}}^2$. We
define the following additional norms:
 \[
\compositenorm{\{\bv,\bV\}}^2=\norm{\bv}_{H^1(\Omega_{\pm})}^2+\norm{\bV}_{H^1(\Gamma)}^2, \quad
 \compositenorm{\{p,\pi\}}^2=\norm{p}^2_{\Omega_{\pm}}+\norm{\pi}^2_{\Gamma}
 \]
The coercivity  result is formulated in the form of a lemma.
\begin{lemma}\label{Lem:coerc} For any $\bv \in V^{\pm}_0$ and $\bV \in V_\Gamma$ it holds
\begin{equation}
a(\{\bv,\bV\},\{\bv,\bV\})\geq C\compositenorm{\{\bv,\bV\}}^2\label{eq:a_coer}
\end{equation}
with a positive constant $C$, which may depend on the viscosity values and $\Omega_\pm$.
\end{lemma}
\begin{proof}
One readily computes that
\begin{multline}\label{aux550}
a(\{\bv,\bV\},\{\bv,\bV\})= 2(\mu^-\bD(\bv^-),\bD(\bv^-))_{\Omega_{-}}  + 2(\mu^+\bD(\bv^+),\bD(\bv^+))_{\Omega_{+}} \\
 +2(\mu_{\Gamma}\bD_{\Gamma}(\bV),\bD_{\Gamma}(\bV))_{\Gamma} + f^-\| \bP\bv^- - \bV \|_{\Gamma}^2  + f^+ \| \bP\bv^+ - \bV \|_{\Gamma}^2.
\end{multline}
Since function $\bv^+$ satisfies homogeneous Dirichlet boundary condition on $\partial\Omega_{+}\setminus\Gamma$,  we apply the following Korn's inequality in $\Omega_{+}$:
\begin{equation}\label{aux555}
\norm{\bv^+}_{H^1(\Omega_{+})}\le C \|\bD(\bv^+)\|_{\Omega_{+}}
\end{equation}
By the triangle and trace inequalities in $\Gamma$, we get
\begin{equation}\label{aux559}
\norm{\bV}_{\Gamma}\le \| \bP\bv^+ - \bV \|_{\Gamma}+ \| \bP\bv^+\|_{\Gamma}\le \| \bP\bv^+ - \bV \|_{\Gamma}+C \norm{\bv^+}_{H^1(\Omega_{+})}.
\end{equation}
We further apply Korn's inequality on $\Gamma$~\cite{jankuhn2017incompressible}:
\begin{equation}\label{aux563}
\norm{\bV}_{H^1(\Gamma)}\le C\left(\norm{\bV}_{\Gamma} + \|\bD_{\Gamma}(\bV)\|_{\Gamma}\right).
\end{equation}
Next, we can estimate the trace of $\bv^-$ on $\Gamma$ through the triangle inequality:
\begin{equation}\label{aux567}
\norm{\bv^-}_{\Gamma} \le \| \bP\bv^- - \bV \|_{\Gamma}+ \| \bV \|_{\Gamma} \le \| \bP\bv^- - \bV \|_{\Gamma}+ \norm{\bV}_{H^1(\Gamma)}.
\end{equation}
We finally apply the following Korn's inequality in $\Omega_{-}$:
\begin{equation}\label{aux571}
\norm{\bv}_{H^1(\Omega_{-})}\le C \left(\|\bD(\bv^-)\|_{\Omega_{-}}+ \norm{\bv^-}_{\Gamma}\right).
\end{equation}
Identity \eqref{aux550} and inequalities \eqref{aux555}--\eqref{aux571} lead to \eqref{eq:a_coer} after easy computations.\\
\end{proof}

The continuity of the bilinear forms $a(\{\cdot,\cdot\},\{\cdot,\cdot\})$, $b(\{\cdot,\cdot\},\{\cdot,\cdot\})$ and $s(\cdot,\cdot)$ follows from standard arguments based on the Cauchy--Schwarz and triangle inequalities:
\begin{equation}\label{cont}
\begin{split}
a(\{\bu,\bU\},\{\bv,\bV\})\le C \compositenorm{\{\bu,\bU\}}\compositenorm{\{\bv,\bV\}}\quad\text{for all}~\bu, \bv \in V^\pm_0,~\bU, \bV \in V_\Gamma, \\
b(\{\bu,\bU\},\{p,\pi\})\le C \compositenorm{\{\bu,\bU\}}\compositenorm{\{p,\pi\}}\quad\text{for all}~\bu\in V^\pm_0,~\bU \in V_\Gamma,~p\in L^2(\Omega),~\pi\in L^2(\Gamma), \\
s(\bv,\pi)\le C \compositenorm{\{\bv,0\}}\compositenorm{\{0,\pi\}}\quad\text{for all}~\bv\in V^\pm,~\pi\in L^2(\Gamma).
\end{split}
\end{equation}
Problem \eqref{weak} falls into the class of so-called generalized saddle point problems.
An abstract well-posedness result for such problems
can be found, e.g.~in~\cite{bernardi1988generalized,nicolaides1982existence}, which
extend the Babu\c{s}ka--Brezzi theory.  Applied to \eqref{weak}, this well-posedness
result requires coercivity~\eqref{eq:a_coer}, continuity~\eqref{cont} and two inf-sup conditions formulated in the following lemma.

\begin{lemma} The following inf-sup conditions hold with positive constants $\gamma_1$ and $\gamma_2$:
\begin{align}
\label{insup1}
    \sup_{\bv\in V^\pm_0, \bV\in V_{\Gamma}}\frac{b(\{\bv,\bV\},\{p,\pi\})+s(\bv,\pi)}{\compositenorm{\{\bv,\bV\}}}\geq\gamma_1\compositenorm{\{p,\pi\}},\quad\forall~p\in L^2_\pm,~\pi\in L^2(\Gamma),  \\
\label{insup2}
     \sup_{\bv\in V^\pm_0, \bV\in V_{\Gamma} }\frac{b(\{\bv,\bV\},\{p,\pi\})}{\compositenorm{\{\bv,\bV\}}}\geq\gamma_2\compositenorm{\{p,\pi\}},\quad\forall~p\in L^2_0(\Omega),~\pi\in L^2_0(\Gamma).
\end{align}
\end{lemma}
\begin{proof} The proof follows by combining well-known results about the existence of a continuous right inverse of the divergence operator in $H^1_0(\Omega)^3$~\cite{bogovskii1979solution} and $V_\Gamma$~\cite{jankuhn2017incompressible}: For arbitrary $p\in L^2_0(\Omega)$ and $\pi\in L^2_0(\Gamma)$ there exist
$\bv\in H^1_0(\Omega)^3$ and $\bV\in V_\Gamma$ such that
\begin{equation}\label{aux610}
\begin{split}
p=\Div \bv~\text{~in}~\Omega,\quad\text{and}\quad \|\bv\|_{H^1(\Omega)}\le c_\Omega\|p\|_{L^2(\Omega)},\\
\pi=\divG \bV~\text{~on}~\Gamma,\quad\text{and}\quad\|\bV\|_{H^1(\Gamma)}\le c_\Gamma\|\pi\|_{\Gamma}.
\end{split}
\end{equation}
Letting $\bv^\pm=\bv|_{\Omega_\pm}$, $(\bv^-,\bv^+)^T \in V^\pm_0$, and adding estimates in \eqref{aux610}  we get
 \begin{align}
\compositenorm{\{p,\pi\}}^2 &\le b(\{\bv,\bV\},\{p,\pi\}), \label{div_stab1a} \\
\compositenorm{\{\bv,\bV \}} & \le c_\Omega\|p\|_{\Omega}+c_\Gamma\|\pi\|_{\Gamma}\le (c_\Omega+c_\Gamma)\compositenorm{\{p,\pi\}}. \label{div_stab2a}
\end{align}
This proves \eqref{insup2} with $\gamma_2=1/(c_\Omega+c_\Gamma)$.
\medskip

To show \eqref{insup1}, we split $\pi=\pi_0+\pi^\perp$ with $\pi_0\in L^2_0(\Gamma)$ and $\pi^\perp=|\Gamma|^{-1}\int_\Gamma \pi\,ds$.
For the $\pi_0$ part of $\pi$, we use again \eqref{aux610} as above, while for $p^\pm\in L^2_0(\Omega^\pm)$ we use the existence of a continuous right inverse of $\Div$ in $H^1_0(\Omega^\pm)^3$ to claim the existence of $\bv\in H^1_0(\Omega^-)^3\times H^1_0(\Omega^+)^3\subset V^\pm_0$ and $\bV\in V_\Gamma$ such that
 \begin{equation}
\compositenorm{\{p,\pi_0\}} \le b(\{\bv,\bV\},\{p,\pi_0\})+s(\bv,\pi), \qquad
\compositenorm{\{\bv,\bV \}} \le (c_\Omega+ c_\Gamma)\compositenorm{\{p,\pi_0\}}, \label{div_stab2}
\end{equation}
with some positive $c_\Omega, c_\Gamma$ depending only on $\Gamma$ and $\Omega$.  We also used that $\bv=\boldsymbol{0}$ on $\Gamma$ implies $s(\bv,\pi)=0$.

Let $C^\pm=\pm|\Omega^\pm|^{-1}|\Gamma|\int_\Gamma \kappa\,ds$. To control $\|\pi^\perp\|_{\Gamma}$, we need $\bv_1\in V^\pm_0$ such that
\begin{equation}\label{aux689}
\Div\bv_1=-C^\pm~~\text{in}~\Omega^\pm,~~ \bv_1\cdot\bn=\kappa~~\text{on}~\Gamma~~\text{and}~~\|\bv_1\|_{H^1(\Omega^\pm)}\le C.
\end{equation}
Such $\bv_1$ can be built, for example, as follows: Let $\bv_1^{-}=\nabla\psi$, where $\psi\in H^2(\Omega^{-})$ solves the Neumann problem
$-\Delta\psi=C^{-}$ in $\Omega^{-}$, $\bn\cdot\nabla\psi=\kappa$ on $\Gamma$. Since $\Gamma=\partial\Omega^{-}$ is smooth, 
by the $H^2$-regularity of the Neumann problem we have that $\|\bv_1^-\|_{H^1(\Omega^-)}\le \|\psi\|_{H^2(\Omega^-)}\le C$.
The boundary $\partial\Omega$ is only Lipschits and so the Neumann problem in $\Omega^+$ is not necessarily $H^2$-regular. To handle this,  we first extend $\bv_1^{-}$ from $\Omega^{-}$ to a function   $\tilde\bv_1$ in $H^1_0(\Omega)^3$ such that $\|\tilde\bv_1\|_{H^1(\Omega^+)}\le c\|\bv_1^-\|_{H^1(\Omega^-)}$ ~\cite{stein1970singular}. Next, we consider
$\bw\in H^1_0(\Omega^+)^3$ such that $\Div\bw=C^+-\Div\tilde\bv_1 \in L^2_0(\Omega^+)$, and $\|\bw\|_{H^1(\Omega^+)}\le c_{\Omega^+}\|\Div\bw\|_{L^2(\Omega^+)}\le C$~\cite{bogovskii1979solution}. The desired $\bv_1^+$ is given in $\Omega^+$ by $\bv_1^+=\tilde\bv_1+\bw$.
Since $\Div\bv_1=-C^\pm$, for $p\in L^2_\pm(\Omega)$ and $\pi\in L^2(\Gamma)$ we have identities
 \begin{equation}\label{aux696}
b(\{\bv_1,0\},\{p,\pi\})=0=(\divG\bV,\pi^\perp)_\Gamma.
 \end{equation}
We also note the equality  $\|\pi^\perp\|_{\Gamma}^2=  \hat c\,   s(\bv_1,\pi^\perp)$, with
\[
\hat c=\pi^\perp |\Gamma|/\int_\Gamma\kappa^2ds.
\]
The denominator above is positive, since $\Gamma$ is closed and so $\kappa$ cannot be zero everywhere on $\Gamma$.
We use \eqref{div_stab1a}--\eqref{aux696} to estimate for some $\beta>0$:
\[
\begin{split}
\compositenorm{\{p,\pi\}}^2&=\compositenorm{\{p,\pi_0\}}^2+\beta\|\pi^\perp\|_{\Gamma}^2\\
&\le b(\{\bv,\bV\},\{p,\pi_0\})+s(\bv,\pi) + s(\beta\hat c\,\bv_1,\pi^\perp)\\
&= b(\{\bv+\beta\hat c\bv_1,\bV\},\{p,\pi\})+s(\bv+\beta\hat c\bv_1,\pi)
-(\beta\hat c\bv_1,\pi_0)_\Gamma\\
&\le b(\{\bv+\beta\hat c\bv_1,\bV\},\{p,\pi\})+s(\bv+\beta\hat c\bv_1,\pi)
+\frac{\beta^2\hat c^2}{2}\|\bv_1\|^2_{\Gamma} +\frac{1}{2}\|\pi_0\|^2_{\Gamma}\\
&\le b(\{\bv+\beta\hat c\bv_1,\bV\},\{p,\pi\})+s(\bv+\beta\hat c\bv_1,\pi)
+c_3\,\beta^2\|\pi^\perp\|_{\Gamma}^2 +\frac{1}{2}\|\pi_0\|^2_{\Gamma}.
\end{split}
\]
with some $c_3>0$ depending only on $\Gamma$ and $\Omega$. For $\beta>0$ sufficiently small such that 
$\frac{\beta}2-c_3\,\beta^2\ge0$, we get
\begin{equation}\label{aux724}
c\compositenorm{\{p,\pi\}}^2\le b(\{\bv+\beta\hat c\bv_1,\bV\},\{p,\pi\})+s(\bv+\beta\hat c\bv_1,\pi),
\end{equation}
with $c>0$ depending only on $\Gamma$ and $\Omega$.
Thanks to the triangle inequality, the second estimate in \eqref{div_stab2} and the definition of $\hat c$ and $\bv_1$, we find the bound
\[
\compositenorm{\bv+\beta\hat c\bv_1,\bV}\le \compositenorm{\bv,\bV}+ \beta\hat c\compositenorm{\beta\hat c\bv_1,0}
\le (c_\Omega+ c_\Gamma)\compositenorm{\{p,\pi_0\}} + C\,\|\pi^\perp\|_{\Gamma} \le C\,\compositenorm{\{p,\pi\}} ,
\]
with $C>0$ depending only on $\Gamma$ and $\Omega$. The combination of the above bound and \eqref{aux724} completes
the proof of the lemma.\qquad
\end{proof}


\subsection{Finite element discretization}\label{sec:FE}

Let $\Omega \subset \R^3$ be a fixed polygonal domain  that strictly contains $\Gamma$.
We consider a family of shape regular tetrahedral triangulations $\{\T_h\}_{h >0}$ of $\Omega$.
We adopt the convention that the elements $T$ and edges $e$ are open
sets and use the over-line symbol to refer to their closure.
Let $h_T=\mbox{diam}(T)$ for $T \in \T_h$.
The set of elements intersecting $\Omega_\pm$ and the set of elements
having a nonzero intersection with $\Gamma$ are
\begin{align}
\T^\pm_h = \{ T \in \T_h:T \cap \Omega_\pm \neq \emptyset \}, \quad
\T_h^\Gamma= \{ T \in \T_h:\overline{T} \cap \Gamma \neq \emptyset \},
\end{align}
respectively.
We assume $\{\T_h^\Gamma\}_{h >0}$ to be quasi-uniform.
The domain formed by all tetrahedra in $\T_h^\Gamma$ is denoted by $\OGamma:=\text{int}(\cup_{T \in \T_h^\Gamma} \overline{T})$.
We define the $h$-dependent domains:
\begin{align}
\Omega_h^\pm = \text{int}\left(\cup_{T \in \T_h^\pm} \overline{T}\right)
\end{align}
and the set of faces of $\T_h^\Gamma$ restricted to the interior of $\Omega_h^\pm$:
\begin{align}
\mathcal{E}_h^{\Gamma,\pm} = \{ e = \text{int}(\partial T_1 \cap \partial T_2):T_1, T_2 \in \T_h^\pm~\text{and}~T_1 \cap \Gamma \neq \emptyset~\text{or}~ T_2 \cap \Gamma \neq \emptyset \}.
\end{align}

For the space discretization of the bulk fluid problems, we restrict our attention to inf-sup stable
finite element pair $\bP_{k+1} - P_k$, $k \geq 1$, i.e. Taylor-Hood elements. Specifically,
we consider the spaces of continuous finite element pressures given by:
\begin{align}
Q_h^- = \{ p \in C(\Omega_h^-): q|_T \in P_k(T)~\forall T \in \T^-_h \}. \el
\end{align}
Space $Q_h^+$ is defined analogously. The trial FE pressure space is given by:
\begin{align}
L^2_\pm(\Omega)_h = \{ p = (p^-, p^+) \in Q_h^- \times Q_h^+\,:\,\int_{\Omega^-}p^-=\int_{\Omega^+}p^+=0 \},\el
\end{align}
and the test space by $Q_h^\pm=Q_h^- \times Q_h^+\cap L^2_0(\Omega)$.
Let
\begin{align}
V_h^- = \{ \bu \in C(\Omega_h^-)^3: \bu|_T \in \bP_{k+1}(T)~\forall T \in \T^-_h \}, \el
\end{align}
with the analogous definition for $V_h^+$. Our FE velocity space is given by:
\begin{align}
V^\pm_h = \{ \bu = (\bu^-, \bu^+) \in (V_h^- \times V_h^+) \}. \el
\end{align}

Functions in $L^2_\pm(\Omega)$ and $V^\pm_h$ and their derivatives are multivalued in $\OGamma$, the overlap of $\Omega_h^-$
and $\Omega_h^+$. The jump of a multivalued function over the
interface is defined as the difference of components coming from $\Omega_h^-$
and $\Omega_h^+$, i.e.~$[ \bu] = \bu^- - \bu^+$ on $\Gamma$. Note that this is the
jump that we have previously denoted with $[ \cdot]^-_+$.
Moreover, we define the following averages:
\begin{equation}
\{ \bu \} = \alpha \bu^+ + \beta \bu^-, \quad
\langle \bu \rangle = \beta \bu^+ + \alpha \bu^-, \label{angle_av}
\end{equation}
where $\alpha$ and $\beta$ are weights to be chosen such that $\alpha+\beta = 1$, $0 \leq \alpha, \beta \leq 1$.
For example, in \cite{CLAUS2019185} the setting $\alpha = \mu_-/(\mu_+ + \mu_-)$ and $\beta = \mu_+/(\mu_+ + \mu_-)$
is suggested. In \cite{caceres2019new}, the authors choose $\alpha=0$, $\beta=1$ if $\mu_{-} \le \mu_{+}$ and $\alpha=1$, $\beta=0$
otherwise.
Below, in \eqref{eq:a_1} and \eqref{eq:b_h} we will use
relationship:
\begin{align}
[ab] = [b] \{a\} + \la b \ra [a]. \label{eq:jump_av}
\end{align}

For the discretization of the surface Stokes problem, we first consider the generalized Taylor--Hood bulk spaces in the strip $\Omega_h^\Gamma$:
\begin{align}
V_{\Gamma,h} &= \{ \bU \in C(\Omega_h^\Gamma)^3: \bU|_T \in \bP_{k+1}(T)~\forall T \in \T^\Gamma_h \}, \cl
Q_{\Gamma,h} &= \{ \pi \in C(\Omega_h^\Gamma): \pi|_T \in P_{k}(T)~\forall T \in \T^\Gamma_h \}, \el
\end{align}
$Q_{\Gamma,h}^0=Q_{\Gamma,h}\cap L^2_0(\Gamma)$.
In the trace finite element method, we use the traces of functions from $V_{\Gamma,h}$ and $Q_{\Gamma,h}$ on $\Gamma$. The inf-sup stability of the resulting trace FEM was analyzed in~\cite{Olshanskii_Reusken2019}  for $k=1$ and extended to higher order  isoparametric trace elements in \cite{jankuhn2020error}. 

In the treatment of the surface Stokes problem, one has to enforce the tangentiality condition $\bU \cdot \bn = 0$ on $\Gamma$. In order to enforce it while avoiding locking, we follow \cite{HANSBO2017298,hansbo2016analysis,jankuhn2017incompressible,reuther2017solving,olshanskii2018finite}
and add a penalty term to the weak formulation. 

A discrete variational analogue of problem \eqref{weak} reads:
Find $(\bu_h, p_h) \in V^\pm_h \times L^2_\pm(\Omega)_h$, and $(\bU_h, \pi_h) \in V_{\Gamma,h} \times Q_{\Gamma,h}$ such that
\begin{equation}\label{FE_formulation}
\left\{
\begin{split}
   a_h(\{\bu_h,\bU_h\},\{\bv_h,\bV_h\})+b_h(\{\bv_h,\bV_h\},\{p_h,\pi_h\})+s_h(\bv_h,\pi_h)&= r_h(\bv_h,\bV_h) \\
   b_h(\{\bu_h,\bU_h\},\{q_h,\tau_h\}) - b_p(p_h, q_h)- b_s(\pi_h,\tau_h) &=0
\end{split}
\right.
\end{equation}
for all $(\bv_h, q_h) \in V_{0,h}^{\pm} \times Q_h^\pm$ and $(\bV_h, \tau_h) \in V_{\Gamma,h} \times Q_{\Gamma,h}^0$.
We define all the bilinear forms in \eqref{FE_formulation} for all $\bu_h \in V_h^{\pm}$, $\bv_h\in V^\pm_{0,h}$, $\bU, \bV \in V_{\Gamma,h}$, $p \in L^2(\Omega)$, $\pi\in L^2(\Gamma)$. Let us start from form $a_h(\{\cdot,\cdot\},\{\cdot,\cdot\})$:
\begin{align}
    a_h(\{\bu_h,\bU_h\},\{\bv_h,\bV_h\})=& a_i(\{\bu_h,\bU_h\},\{\bv_h,\bV_h\}) + a_n(\bu_h,\bv_h) \cl
    &+ a_p(\{\bu_h,\bU_h\},\{\bv_h,\bV_h\}) +a_s(\bU_h,\bV_h), \label{eq:a_h}
\end{align}
where we group together the terms that arise from the integration by parts of the divergence of the stress tensors:
\begin{align}
    a_i(\{\bu_h,\bU_h\},\{\bv_h,&\bV_h\})= 2(\mu^-\bD(\bu_h^-),\bD(\bv_h^-))_{\Omega_{-}}  + 2(\mu^+\bD(\bu_h^+),\bD(\bv_h^+))_{\Omega_{+}} \cl
    & + \langle f^-(\bP\bu_h^--\bU_h ), \bP\bv_h^- - \bV_h \rangle_\Gamma
     + \langle f^+(\bP\bu_h^+-\bU_h ), \bP\bv_h^+ - \bV_h  \rangle_\Gamma \cl
    &- 2 \langle \{ \mu \bn^T \bD(\bu_h) \bn \}, [ \bv_h \cdot \bn] \rangle_\Gamma +2(\mu_{\Gamma}\bD_{\Gamma}(\bU_h),\bD_{\Gamma}(\bV_h))_{\Gamma}   \label{eq:a_1}
\end{align}
the terms that enforce condition \eqref{eq:scc1} weakly using Nitsche's method
\begin{align}
a_n(\bu_h,\bv_h)= \sum_{T \in \mathcal{T}_h^\Gamma} \frac{\gamma}{h_T} \{ \mu \} ([\bu_h \cdot \bn], [\bv_h \cdot \bn])_\Gamma -
2 \langle \{ \mu \bn^T \bD(\bv) \bn \} ,  [\bu_h \cdot \bn] \rangle_\Gamma, \label{eq:a_n}
\end{align}
and the stabilization and penalty terms:
\begin{align}
a_p(\{\bu_h,\bU_h\},\{\bv_h,\bV_h\}) &= \bJ_h^-(\bu_h, \bv_h) + \bJ_h^+(\bu_h, \bv_h) + \tau_s (\bU_h \cdot \bn, \bV_h \cdot \bn)_\Gamma, \label{eq:a_p} \\
\bJ_h^\pm(\bu_h, \bv_h) &= \sum_{\ell = 1}^{k+1} | e|^{2\ell - 1}  \sum_{e \in \mathcal{E}_h^{\Gamma, \pm}} \gamma_\bu^\pm \mu^\pm ([\partial_n^\ell \bu_h^\pm], [\partial_n^\ell \bu_h^\pm])_e. \label{eq:a_j}
\end{align}
In \eqref{eq:a_j}, $\partial_n^\ell \bu_h^-$ denotes the derivative of order $\ell$ of $\bu_h^-$ in the direction of $\bn$.
The $\bJ_h$ terms in \eqref{eq:a_p} are so called ghost-penalty stabilization~\cite{B10,cutFEM}  included to avoid poorly conditioned algebraic systems due  to possible small cuts of tetrahedra from $\T^\Gamma_h$ by the interface.
The terms in \eqref{eq:a_s} and \eqref{eq:b_s} have the same role  for the surface bilinear forms.

The last form in \eqref{eq:a_h} is related to the algebraic stability of the surface Stokes problem:
\begin{align}
a_s(\bU_h,\bV_h) = \rho_\bu (\nabla \bu_h \bn, \nabla \bv_h \bn)_{\Omega_h^\Gamma}. \label{eq:a_s}
\end{align}

Similarly, the terms coming from the integration by parts of the divergence of the stress tensors are contained in
\begin{align}
        b_h(\{\bv_h,\bV_h\},\{p_h,\pi_h\})=&- (p_h^-,\Div  \bv_h^-)_{\Omega_{-}} - (p_h^+,\Div  \bv_h^+)_{\Omega_{+}} \cl
    &+ \la \{ p_h \}, [ \bv_h \cdot \bn] \ra_\Gamma +(\nablaG\pi_h,\bV_h)_{\Gamma}, \label{eq:b_h}
\end{align}
the penalty terms are grouped together in
\[
b_p(p_h, q_h) = J_h^-(p_h, q_h) + J_h^+(p_h, q_h), \quad
J_h^\pm(p_h, q_h) = \frac{\gamma_p^\pm}{\mu^\pm} \sum_{e \in \mathcal{E}_h^{\Gamma, \pm}} \sum_{\ell = 1}^k |e|^{2\ell + 1} ([\partial_n^\ell p_h^\pm], [\partial_n^\ell q_h^\pm])_e,
\]
and we have a term related to algebraic stability of the surface Stokes problem in:
\begin{align}
b_s(\pi_h,\tau_h) = \rho_p (\nabla p_h \cdot \bn, \nabla p_h \cdot \bn)_{\Omega_h^\Gamma}. \label{eq:b_s}
\end{align}
Finally,
\begin{align}
     s_h(\bv_h,\pi_h)=&-\langle \pi_h \kappa , \la \bv_h \cdot \bn  \ra \rangle_\Gamma,\cl
     r_h(\bv_h,\bV_h)=&(\bbf_h^-,\bv_h^-)_{\Omega_{-}} + (\bbf_h^+,\bv_h^+)_{\Omega_{+}}+(\bP\bb_h^e,\bV_h)_{\Gamma}. \el
\end{align}
We recall that some of the interface terms in $a_i(\{\cdot,\cdot\},\{\cdot,\cdot\})$ and $b_h(\{\cdot,\cdot\},\{\cdot,\cdot\})$
have been obtained using relationship \eqref{eq:jump_av}.

Parameters $\gamma_\bu^\pm$, $\gamma_p^\pm$, and $\gamma$ are all assumed to be independent
of $\mu^\pm$, $h$, and the position of $\Gamma$ against the underlying mesh. Parameter $\gamma$
in \eqref{eq:a_n} needs to be large enough to provide the bilinear form $a_h(\{\cdot,\cdot\},\{\cdot,\cdot\})$
with coercivity. Parameters  $\gamma_\bu^\pm$ and $\gamma_p^\pm$ can be tuned to improve
the numerical performance of the method.
As for the parameters required by the discretization of the surface Stokes problem, we allow:
\begin{align}
\tau_s = c_\tau h^{-2}, \quad \rho_p = c_p h, \quad \rho_\bu \in [c_\bu h, C_\bu h^{-1}],
\end{align}
where $c_\tau$, $c_p$, $c_\bu$, and $C_\bu$ are positive constants independent of
$h$ and how $\Gamma$ cuts the bulk mesh.

The definition of bilinear forms requires integration over $\Gamma\cap T$ and  $T\cap\Omega^\pm$ for $T$ from $\Omega_h^\Gamma$. In general, there are no exact quadrature formulas to accomplish this task~\cite{OlshSafin2}. In practice, approximations should be made which introduce geometric errors. To keep these geometric errors of the order consistent with the approximation properties of the finite element spaces,  we use isoparametric  variants of the above spaces introduced in~\cite{lehrenfeld2016high}; see also \cite{lehrenfeld2018analysis,grande2018analysis}.

We expect that the stability of the finite element formulation can be analyzed largely following 
the same steps of the well-posedness analysis for the weak formulation in Sec.~\ref{sec:well}, 
with a special treatment of cut elements, Nitsche terms and surface elements as available in the literature for bulk 
Stokes interface and surface Stokes problems. For the sake of brevity,  we do not work out these details here, but will 
present them in a follow-up paper.

\section{A partitioned method for the coupled bulk-surface flow}\label{sec:part_meth}

For the solution of the coupled problem described in Sec.~\ref{sec:simplified_p} we intend to use a partitioned strategy, i.e.~each
sub-problem is solved separately and the coupling conditions are enforced in an iterative fashion.
Partitioned method are appealing for solving coupled problems because they allow to reuse existing solvers with minimal modifications.
In order to devise such a method for the simplified problem in Sec.~\ref{sec:simplified_p},
let us take a step back and look at the original problem \eqref{eq:Full}.

Discretize  problem \eqref{eq:Full} in time with, e.g., the Backward Euler method  and consider the coupled problem at
a particular time $t = t^{n+1}$. Let $S_b$ be the map that associates the jump in the normal stress across
the interface to any given surface flow velocity $\bU=\bU_T+U_N\bn$:
\[
S_b(\bU) = [\bsigma \bn]^+_{-} = \bsigma^+(\bu^+, p^+) \bn - \bsigma^-(\bu^-, p^-) \bn \quad \text{on}~\Gamma,
\]
where $(\bu^+, p^+)$ and $(\bu^-, p^-)$ represent the solution of the two-phase time-discrete Navier-Stokes problem
at time $t$  associated to \eqref{eq:NS1}--\eqref{eq:NS2} endowed with interface conditions \eqref{eq:cc1}, \eqref{eq:cc2},
and \eqref{eq:cc3}. Moreover, let $S_s$ be the operator associated to the surface flow such that to any given
surface flow velocity $\bU$ it associates the load $\bbf_\Gamma$:
\[
S_s(\bU) = \bbf_\Gamma \quad \text{on}~\Gamma,
\]
through the time-discrete surface Navier-Stokes problem
at time $t$ associated to \eqref{eq:sNS1}--\eqref{eq:sNS2}. Note that $S_b$ and $S_s$ are nonlinear and
their definitions can involve also forcing terms and, in the case of the bulk fluid problem, terms due to the boundary conditions.
For the surface operator, we can define $S^{-1}_s$ as the map that associates the surface flow velocity $\bU$
to any given load $\bbf_\Gamma$ on $\Gamma$.

With the above definitions, we can express the time discrete version of coupled problem \eqref{eq:Full}
in terms of the solution $\bU$ of a nonlinear equation defined only on $\Gamma$. This interface equation is usually
presented in one of three formulations that are equivalent from the mathematical point of view, but give rise to different
iterative algorithms. The first and perhaps most used formulation is the fixed-point one: Find $\bU$ such that
\begin{equation}\label{eq:fixed-point}
S_s^{-1}(S_b(\bU)) = \bU \quad \text{on}~\Gamma.
\end{equation}
The second formulation is a slight modification of \eqref{eq:fixed-point}, which lends itself to
a Newton iterative method: Find $\bU$ such that
\begin{equation*}
S_s^{-1}(S_b(\bU)) - \bU = \boldsymbol{0} \quad \text{on}~\Gamma.
\end{equation*}
The third approach is given by the Steklov-Poincar\'{e} equation: Find $\bU$ such that
\begin{equation*}
S_b(\bU) - S_s(\bU) = \boldsymbol{0} \quad \text{on}~\Gamma.
\end{equation*}
See, e.g., \cite{B-quarteroniv1} for more details on these three formulations.

A standard algorithm for eq.~\eqref{eq:fixed-point} uses fixed-point iterations: Given $\bU^k$, compute
\begin{equation}\label{alg:fixed-point}
\bU^{k+1} = \bU^{k} + \omega^k (\overline{\bU}^{k} - \bU^{k}) \quad \text{with}~\overline{\bU}^{k} =  S_s^{-1}(S_b(\bU^k)).
\end{equation}
The choice of the relaxation parameter $\omega^k$ determines the efficiency of the algorithm or it might be crucial for convergence
in certain ranges of the physical parameters. An effective strategy for setting $\omega^k$ is the Aitken's acceleration method.

For simplicity, we present algorithm \eqref{alg:fixed-point} applied to the time discrete version of coupled problem \eqref{eq:Full}
with $\omega^k = 1$ for all $k$ (i.e., no relaxation). At time $t = t^{n+1}$, assuming that $\bU ^k$ is known, perform the following steps:
\begin{itemize}
\item[-] {\bf Step 1}: solve the two-phase time-discrete Navier-Stokes problem
at time $t$  associated to \eqref{eq:NS1}--\eqref{eq:NS2} for the bulk flow variables $(\bu^-_{k+1}, p^-_{k+1})$ and $(\bu^+_{k+1}, p^+_{k+1})$
with interface conditions
\begin{align}
\bu^+_{k+1}\cdot \bn&= U_N^k =  \bu^-_{k+1}\cdot \bn  &&\text{on }\Gamma \cl
\bP{\bsigma^+_{k+1}\bn}&=f^+(\bP\bu^+_{k+1}-\bU_T ^k)  &&\text{on }\Gamma,  \cl
\bP{\bsigma^-_{k+1}\bn}&=-f^-(\bP\bu^-_{k+1}-\bU_T^k)  &&\text{on }\Gamma.  \el
\end{align}
\item[-] {\bf Step 2}: solve the time-discrete surface Navier-Stokes problem
at time $t$ associated to \eqref{eq:sNS1}--\eqref{eq:sNS2} for variables $(\bU^{k+1}, \pi^{k+1})$
with interface condition
\begin{align}
\bbf_\Gamma^{k+1} &= [\bsigma_{k+1}\bn]^+_- &\text{ on }\Gamma.  \el
\end{align}
\item[-] {\bf Step 3}: Check the stopping criterion
\begin{align}
|| \bU ^{k+1} - \bU ^{k} ||_{\Gamma} < \epsilon|| \bU ^{k} ||_{\Gamma}, \el
\end{align}
where $\epsilon$ is a given stopping tolerance.
\end{itemize}
Notice that the bulk and surface flow problems are solved separately and sequentially. In general, this algorithm is easy
to implement but convergence could be slow in certain ranges of the physical parameters and require relaxation for speed-up.

The above algorithm adapted to simplified problem \eqref{eq:Stokes1}--\eqref{eq:scc4a},
reads as follows. At iteration $k+1$, assuming that $(\bU_T^k,\pi^k)$ are known, perform the following steps:
\begin{itemize}
\item[-] {\bf Step 1}: solve two-phase problem \eqref{eq:Stokes1}--\eqref{eq:Stokes2} for the bulk flow variables $(\bu^-_{k+1}, p^-_{k+1})$ and $(\bu^+_{k+1}, p^+_{k+1})$
with interface conditions
\begin{align}
\bu^+_{k+1}\cdot \bn&= \bu^-_{k+1}\cdot \bn  &&\text{on }\Gamma \label{sic_1} \\
\bP{\bsigma^+_{k+1}\bn}&=f^+(\bP\bu^+_{k+1}-\bU_T^k)  &&\text{on }\Gamma,  \label{sic_2} \\
\bP{\bsigma^-_{k+1}\bn}&=-f^-(\bP\bu^-_{k+1}-\bU_T^k)  &&\text{on }\Gamma,  \label{sic_3} \\
\left[ \bn^T \bsigma_{k+1}\bn \right]^-_+&=\pi^k\kappa&&\text{on }\Gamma. \label{sic_4}
\end{align}
\item[-] {\bf Step 2}: solve surface flow problem \eqref{eq:sStokes1}--\eqref{eq:sStokes2} for variables $(\bU_T^{k+1}, \pi^{k+1})$
with interface condition
\begin{align}
\bP\bbf_\Gamma^{k+1} &= [\bP\bsigma_{k+1}\bn]^+_- &\text{ on }\Gamma.  \label{sic_5}
\end{align}
\item[-] {\bf Step 3}: Check the stopping criterion
\begin{align}
|| \bU_T^{k+1} - \bU_T^{k} ||_{\Gamma} < \epsilon|| \bU_T^{k} ||_{\Gamma}. \label{eq:stop}
\end{align}
\end{itemize}

Notice that only interface condition \eqref{sic_2}--\eqref{sic_5} are coupling conditions for bulk and surface flows.
If one was to compute the load exerted on the surface fluid in \eqref{sic_5} directly from the solution of the
problem at Step 1, the overall accuracy of the method would be spoiled. Instead, one can compute
$\bP\bbf_\Gamma^{k+1}$ by plugging \eqref{sic_2}--\eqref{sic_3} into \eqref{sic_5}:
\begin{align}
\bP\bbf_\Gamma^{k+1} &= f^+\bP\bu^+_{k+1} + f^-\bP\bu^-_{k+1} - (f^+ + f^-)\bU_T^k   &\text{ on }\Gamma. \el 
\end{align}
However, we prefer to use a more implicit version of the above condition:
\begin{align}
\bP\bbf_\Gamma^{k+1} &= f^+\bP\bu^+_{k+1} + f^-\bP\bu^-_{k+1} - (f^+ + f^-)\bU_T^{k+1}   &\text{ on }\Gamma. \el  
\end{align}
since it could help have a better control of approximate rigid rotations (Killing vector fields).


\section{Numerical results}\label{sec:num_res}

The aim of the numerical results collected in this section is to provide evidence of the robustness of the proposed finite element approach
with respect to the contrast in viscosity in the bulk fluid, surface fluid viscosity, value of the slip coefficients, and position of the interface relative to the fixed computational mesh.

For the averages in \eqref{angle_av}, we set $\alpha=0$ and $\beta=1$ for all the numerical experiments
since we have $\mu_{-}\le\mu_+$.
In addition, we set $\gamma^\pm_\bu = 0.05$, $\gamma_p^\pm = 0.05$, and $\gamma=80$.
The value of all other parameters will depend on the specific test. The stopping tolerance for
criterion \eqref{eq:stop} is set to $\epsilon = 10^{-6}$.
For all the simulations, we choose to use finite element pair $\bP_{2} - P_1$ for both the bulk and surface fluid
problems.

For all the results presented below, we will report the $L^2$ error and a weighted $H^1$ error for the bulk velocity defined as
\begin{equation}\label{eq:weighted_H1}
\left(2\mu_{-}\|D(\bu-\bu_h^-)\|_{\Omega^-}^2+2\mu_{+}\|D(\bu-\bu_h^+)\|_{\Omega^+}^2\right)^{\frac12},
\end{equation}
and a weighted $L^2$ error for the bulk pressure defined as
\begin{equation}\label{eq:weighted_L2}
\left(
\mu^{-1}_{-}\|p-p_h^-\|_{\Omega^-}^2 + \mu^{-1}_{+}\|p-p_h^+\|_{\Omega^+}^2
\right)^{\frac12}.
\end{equation}
Such weighted norm naturally arise in the error analysis of the Stokes interface problem~\cite{Olshanskii04}.
In addition, we will report the  $L^2$ and $H^1$ errors for the surface velocity and  $L^2$ error for the surface pressure.

\subsection{Sphere embedded in a cube}

We perform a series of tests where domain $\Omega$ is the cube $[-1.5,1.5]^3$ and interface $\Gamma$
is the unit sphere centered at the origin. Let $\bx = (x,y,z) \in \Omega$.
Surface $\Gamma$ is characterized as the zero level set of function $\phi(\bx)= || \bx ||^2_2 - 1$.
We consider the following  solution for the bulk
flow:
\begin{align}
p^-&=3x\sqrt{x^2+y^2+z^2}-2x(x^2+y^2+z^2),  \quad \bu^-=\frac{2f^-}{f^--\mu^-} \ba(x,y,z), \label{ex_p_sphere} \\
p^+&=6x\sqrt{x^2+y^2+z^2}-4x(x^2+y^2+z^2) , \quad \bu^+ = \frac{2f^+}{f^++\mu^+}\ba(x,y,z), \label{ex_u_sphere}
\end{align}
where
\begin{align}
\ba(x,y,z) = \left(\frac{3}{2}-\sqrt{x^2+y^2+z^2}\right)
\left[
\begin{array}{c}
(-y-z)x+y^2+z^2  \\
(-x-z)y+x^2+z^2 \\
(-x-y)z+y^2+x^2
\end{array}
\right], \el
\end{align}
coupled to the following exact solution for the surface flow:
\begin{align}
\pi=x, \quad
\bU  =
\left[
\begin{array}{c}
(-y-z)x+y^2+z^2  \\
(-x-z)y+x^2+z^2 \\
(-x-y)z+y^2+x^2
\end{array}
\right], \label{ex_surface_sphere}
\end{align}
The forcing terms $\bbf^-$ and $\bbf^+$ are found by plugging the solution \eqref{ex_p_sphere}--\eqref{ex_u_sphere} in \eqref{eq:Stokes1}.
We impose a Dirichlet condition \eqref{eq:bcD} on the faces $x=1.5$, $y=-1.5$, $z=-1.5$, where function $\bg$
is found from $\bu^+$ in \eqref{ex_u_sphere}. On the remaining part of the boundary, we impose
a Neumann condition \eqref{eq:bcN} where $\blf_N$ is found from $p^+$ in \eqref{ex_p_sphere} and $\bu^+$ in \eqref{ex_u_sphere}.

The value of the physical parameters will be specified for each test.

\begin{figure}[hbt!]
\centering
\hskip -.3cm
\includegraphics[width=.56\textwidth]{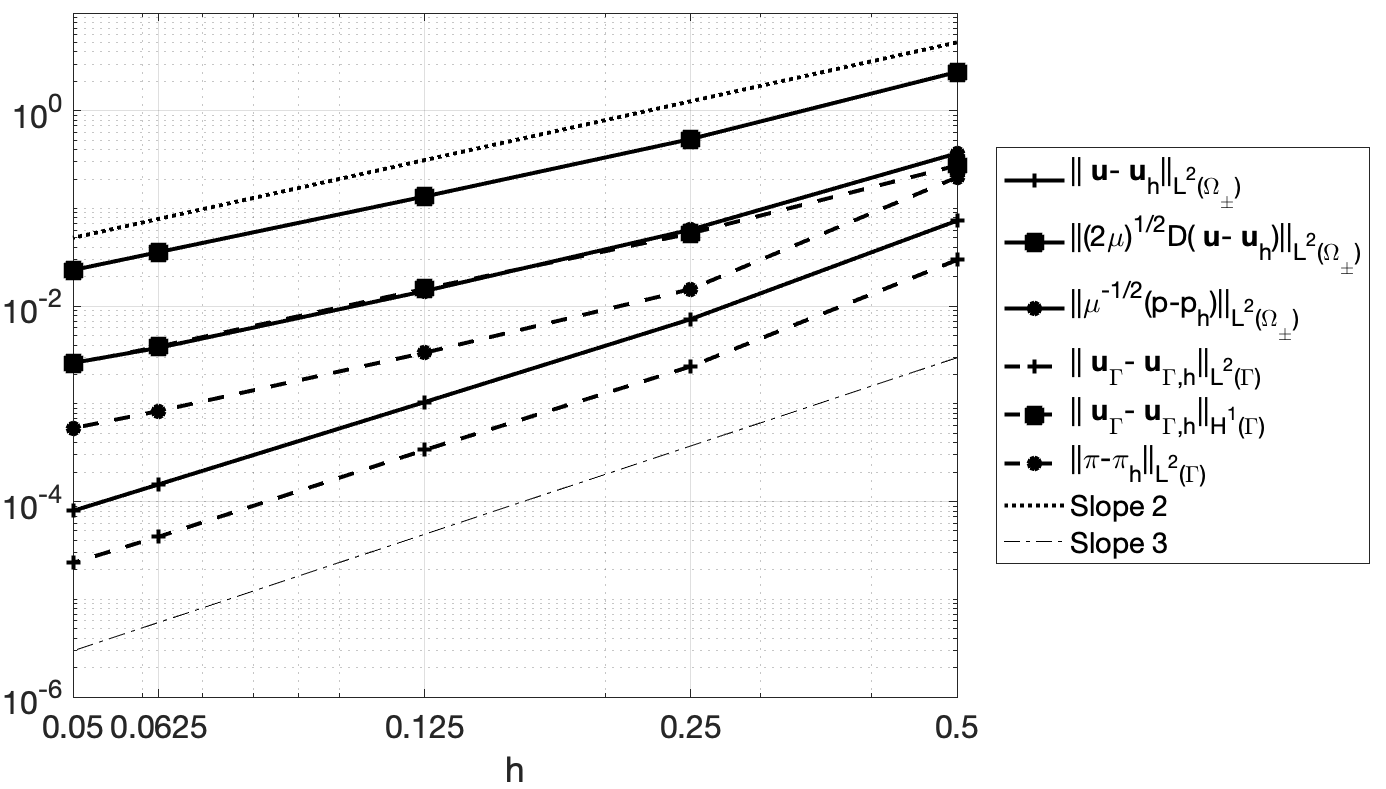}\quad
\includegraphics[width=.39\textwidth]{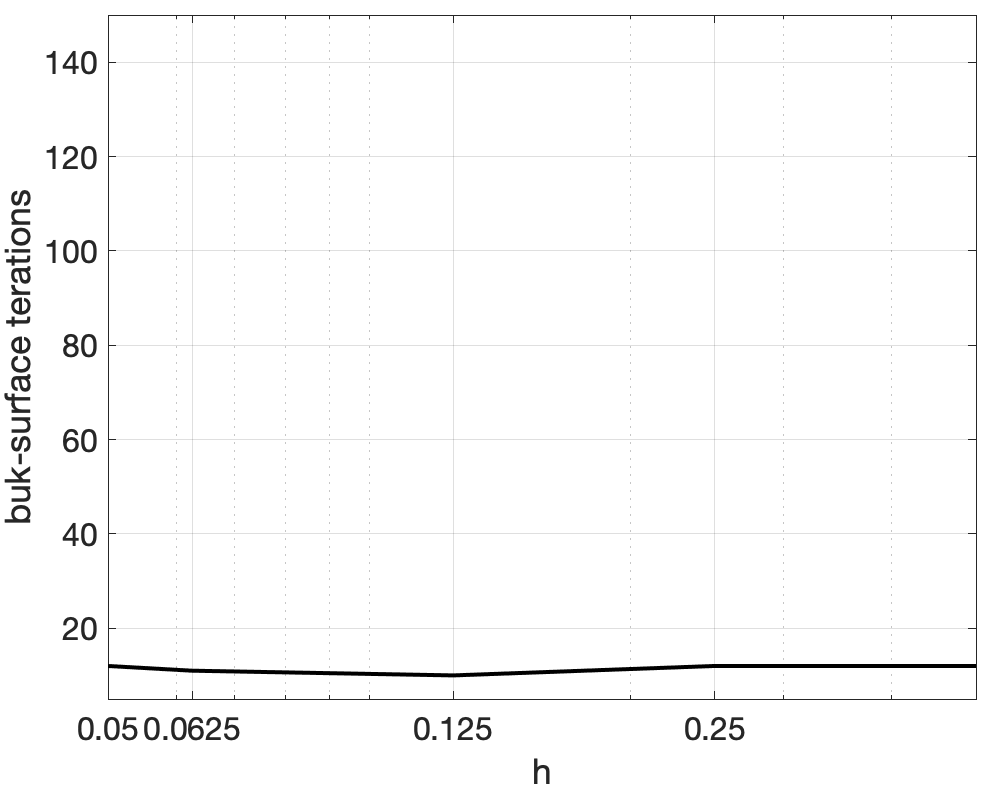}
  \caption{Sphere: (left) Bulk and surface FE errors against the mesh size $h$. (right) Number of bulk-surface iterations of the partitioned method as $h$ varies.}
  \label{fig:spatial_sphere}
\end{figure}

\vskip .2cm
\noindent {\bf Spatial convergence.}
To check the spatial accuracy of the finite element method described in Sec.~\ref{sec:FE},
we consider exact solution \eqref{ex_p_sphere}--\eqref{ex_surface_sphere} with  viscosities $\mu^-=1$, $\mu^+=10$
and $\mu_\Gamma = 1$, and friction coefficients $f^-=2$ and $f^+=10$. Notice that the fluid outside the sphere has
a larger viscosity than the fluid inside the sphere, which has the same viscosity as the surface fluid.
We consider structured meshes of tetrahedra with five levels of refinement, the coarsest mesh
having mesh size $h = 0.5$ while the finest mesh has $h = 0.05$.
All the meshes feature a local one-level refinement near the corners of $\Omega$.
Table \ref{tab:meshes_sphere} reports the number of DOFs for each mesh.
Fig.~\ref{fig:spatial_sphere} (left) shows the $L^2$ error and weighted $H^1$ error
\eqref{eq:weighted_H1} for the bulk velocity, weighted $L^2$ error  \eqref{eq:weighted_L2}
for the bulk pressure,  $L^2$ and $H^1$ errors for the surface velocity and  $L^2$ error for the surface pressure
against the mesh size $h$.
We observe optimal convergence rates for all the norms under consideration.
Fig.~\ref{fig:spatial_sphere} (right) shows the number of bulk-surface iterations to satisfy stopping
criterion \eqref{eq:stop} as $h$ varies. As we can see, the number of iterations
is fairly insensitive to a mesh refinement or coarsening.

\begin{table}[h]
\centering
\begin{tabular}{|c|c|c|c|c|c|}
\hline
$h$ & 0.5 & 0.25 & 0.125 & 0.0625 & 0.05  \\
\hline
\# bulk velocity DOFs & $1.1e4$ & $7.4e4$ & $5.2e5$ & $3.6e6$ & $6.4e6$  \\
\hline
\# bulk pressure DOFs & $6.2e2$ & $3.7e3$ & $2.3e4$ & $1.6e5$ & $2.8e5$  \\
\hline
\# surface velocity DOFs & $2.4e3$ & $1.0e4$ & $4.0e4$ & $1.5e5$ & $2.2e5$  \\
\hline
\# surface pressure DOFs & $1.4e2$ & $5.9e2$ & $2.3e3$ & $8.5e3$ & $1.3e4$  \\
\hline
\end{tabular}
\caption{Sphere: DOFs for bulk and surface variables for all the meshes under consideration in the spatial convergence test.}
\label{tab:meshes_sphere}
\end{table}

{\bf Robustness with respect to the viscosity contrast.} It is known that the case
of high contrast for the viscosities in a two-phase problem is especially challenging from the
numerical point of view. To test the robustness of our approach with respect to the viscosity contrast in the bulk,
we consider exact solution \eqref{ex_p_sphere}--\eqref{ex_surface_sphere} and fix $\mu^-=1$, while we let
$\mu^+$ vary from 1 to 256. We set $\mu_\Gamma = 1$ and friction coefficients $f^-=2$ and $f^+=10$.

\begin{figure}[hbt!]
\centering
\hskip -.3cm
\includegraphics[width=.56\textwidth]{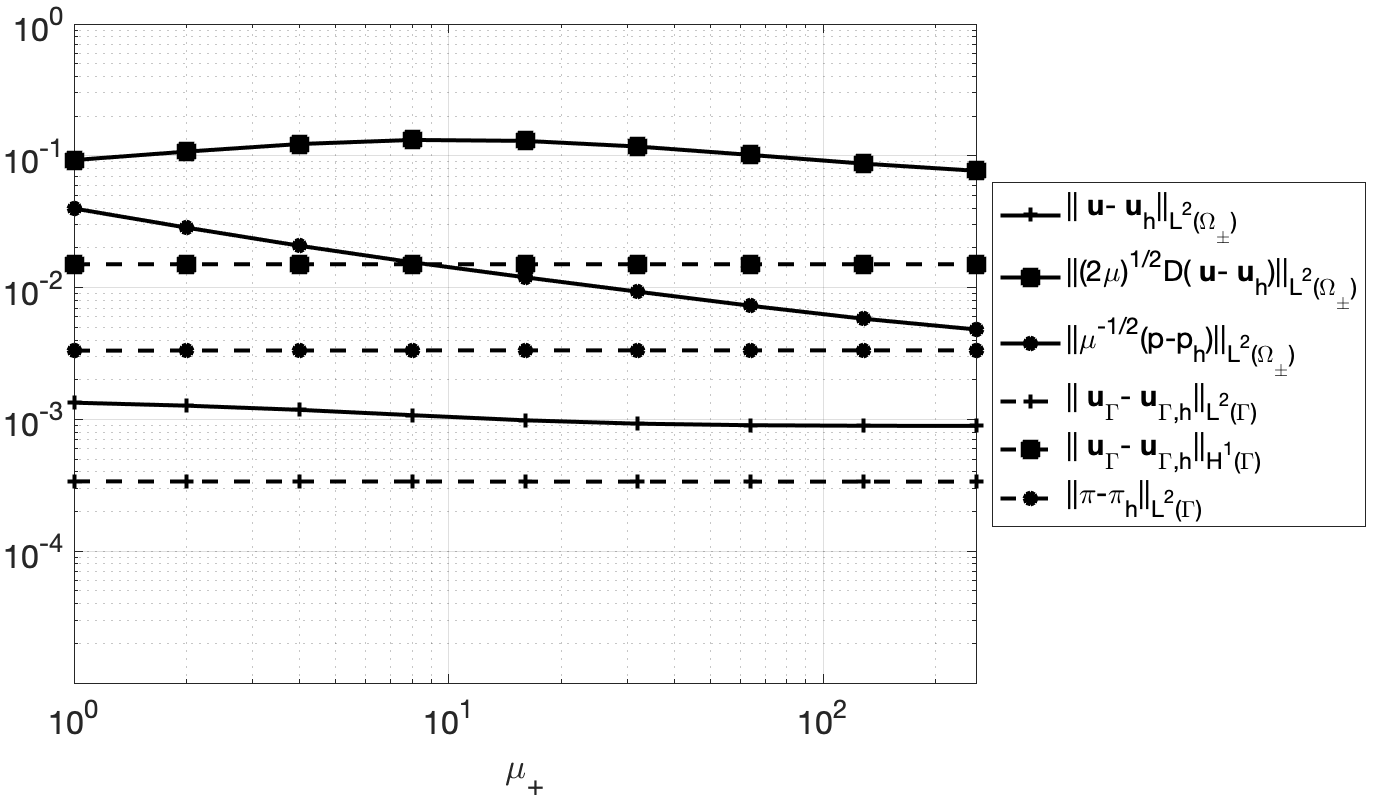}\quad
\includegraphics[width=.39\textwidth]{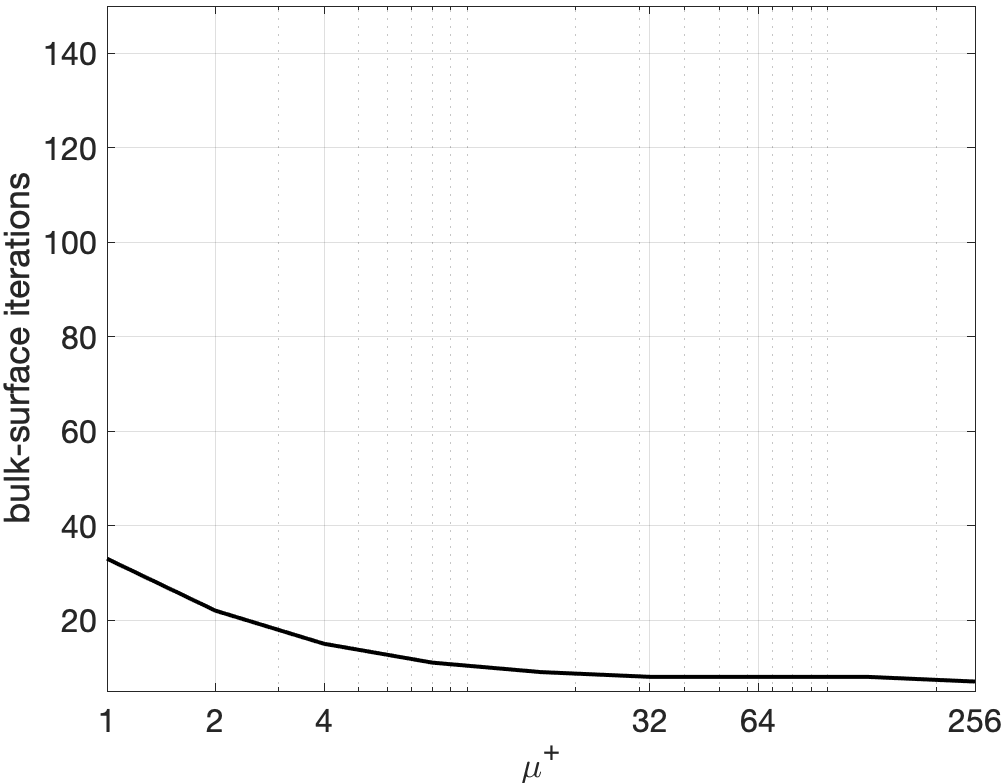}
  \caption{Sphere: (left) Bulk and surface FE errors against the value of $\mu^+$. (right) Number of bulk-surface iterations of the partitioned method as $\mu^+$ varies.}
  \label{fig:muplus_sphere}
\end{figure}

We consider one of the meshes adopted for the previous sets of simulations (with $h = 0.125$).
Fig.~\ref{fig:muplus_sphere} (left) shows the $L^2$ error and weighted $H^1$ error
\eqref{eq:weighted_H1} for the bulk velocity, weighted $L^2$ error  \eqref{eq:weighted_L2}
for the bulk pressure,  $L^2$ and $H^1$ errors for the surface velocity and  $L^2$ error for the surface pressure
against the value of $\mu^+$. We see that the errors remain mostly unchanged as $\mu^+$ varies, with the exception of
the weighted $L^2$ error for the bulk pressure, which decreases as $\mu^+$ increases.
In \cite{olshanskii2021unfitted}, which focuses only on two-phase bulk flow, we found that such error reaches a plateau
as  $\mu^+$ is further increased.
Fig.~\ref{fig:muplus_sphere} (left) shows that our approach is substantially robust with respect to the
viscosity contrast $\mu^+/\mu^-$.

Fig.~\ref{fig:muplus_sphere} (right) reports the number of bulk-surface iterations to satisfy stopping
criterion \eqref{eq:stop} as $\mu^+$ varies. We observe that the number of iterations increases
as the $\mu^+/\mu^-$ ratio decreases, indicating that the coupled bulk-surface problem becomes more
stiff as $\mu^+$ decreases to match $\mu^-$ and $\mu_\Gamma$.

\begin{figure}[hbt!]
\centering
\hskip -.3cm
\includegraphics[width=.56\textwidth]{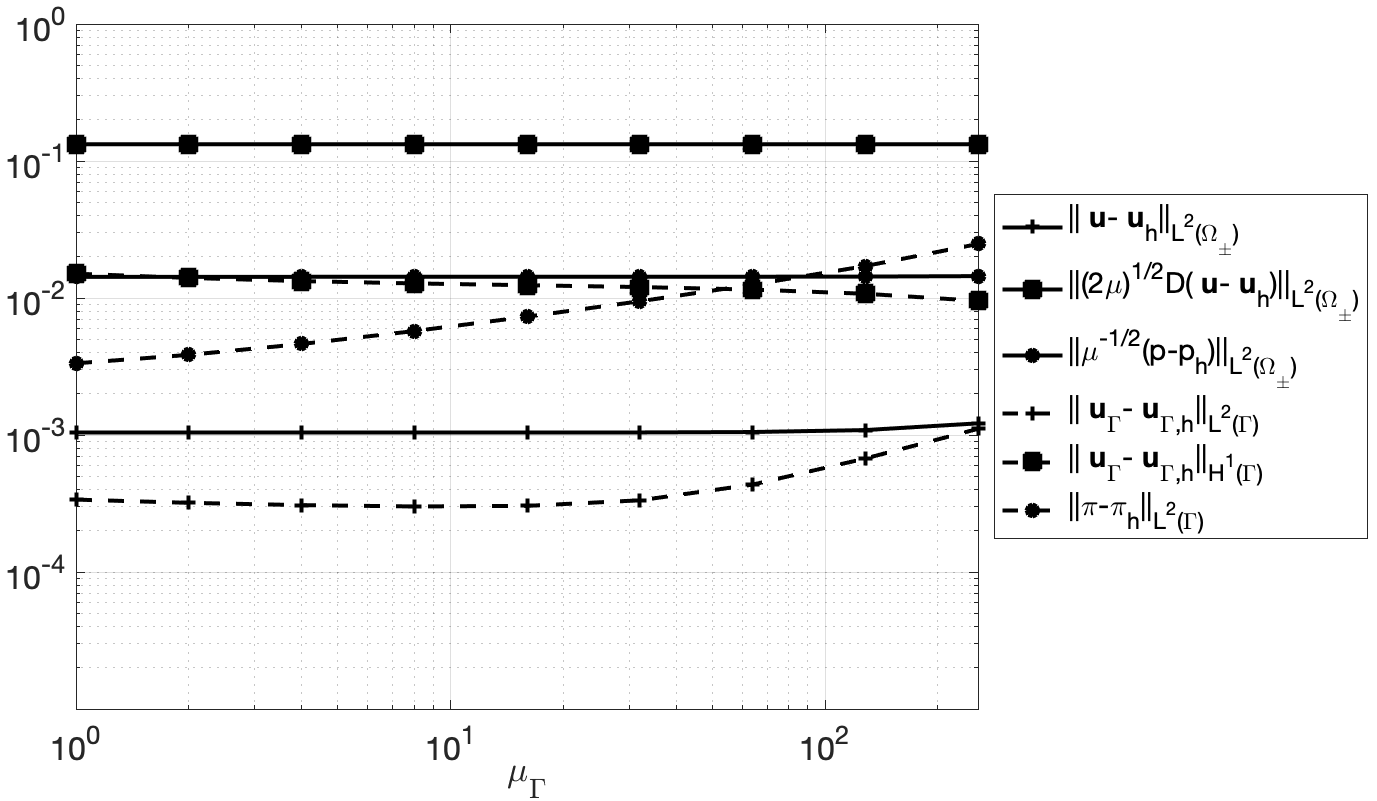}\quad
\includegraphics[width=.39\textwidth]{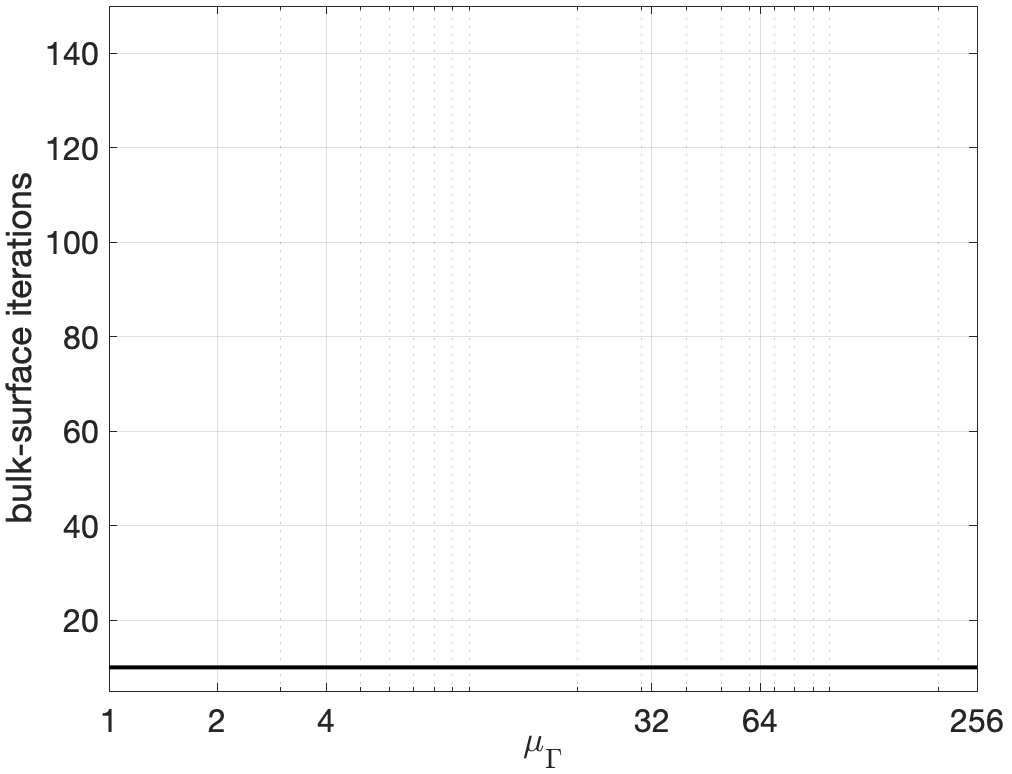}
  \caption{Sphere: (left) Bulk and surface FE errors against the value of $\mu_\Gamma$. (right) Number of bulk-surface iterations of the partitioned method as $\mu_\Gamma$ varies.}
  \label{fig:muGamma_sphere}
\end{figure}

{\bf Robustness with respect to the value of the surface viscosity.} We now let $\mu_\Gamma$ vary from 1 to 256 and keep
all the other physical parameters fixed to the following values:  $\mu^-=1$, $\mu^+=10$, $f^-=2$ and $f^+=10$. Again, we consider exact solution \eqref{ex_p_sphere}--\eqref{ex_surface_sphere} and the mesh with mesh size $h = 0.125$.
Fig.~\ref{fig:muGamma_sphere} (left) shows all the errors we have considered
so far against the value of $\mu_\Gamma$. We notice that all the bulk errors stay constant as $\mu_\Gamma$ varies. The
$L^2$ errors for the surface velocity and pressure increase as $\mu_\Gamma$ increases, while the $H^1$ error for the
surface velocity slightly decreases as $\mu_\Gamma$ increases. This experiment suggests that more viscous embedded layer is less controlled by the bulk fluid which effects the numerical stability of the complete system. In a water -- lipid membrane system, the ratio of lateral dynamic viscosities of the embedded bi-layer and bulk water is typically  1--10 $\mu m$  (depending on the temperature and composition) with the size of a vesicle being generally between 0.1 and 10 $\mu m$. Hence the observed  increase of the numerical error does not look critical for this application.

Fig.~\ref{fig:muGamma_sphere} (right) shows the number of bulk-surface iterations to satisfy stopping
criterion \eqref{eq:stop} as $\mu_\Gamma$ varies. Our partitioned method seems to be insensitive to a
variation in the value of $\mu_\Gamma$. In particular, for the range of  $\mu_\Gamma$ under consideration
the number of iterations stays constant and equal to 12.

\begin{figure}[hbt!]
\centering
\hskip -.3cm
\includegraphics[width=.56\textwidth]{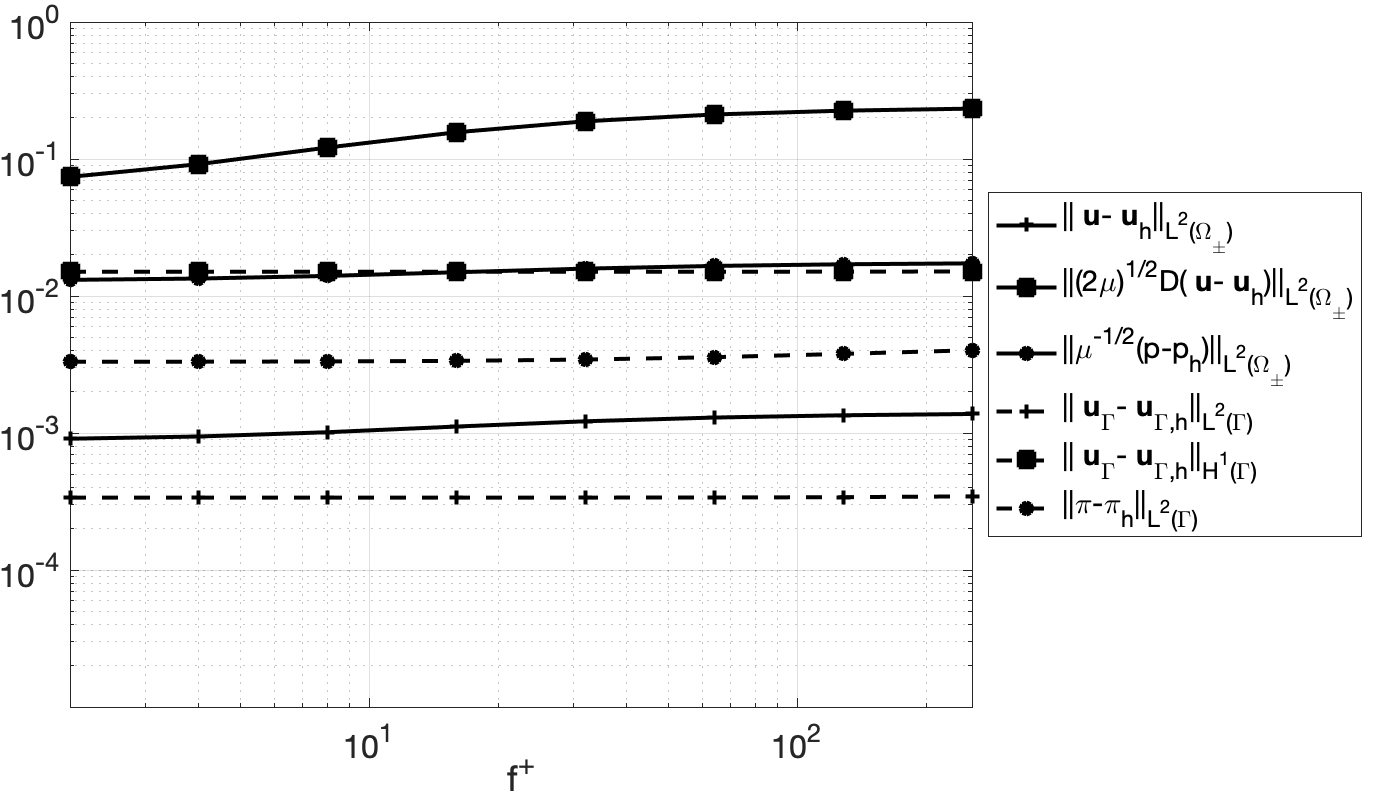}\quad
\includegraphics[width=.39\textwidth]{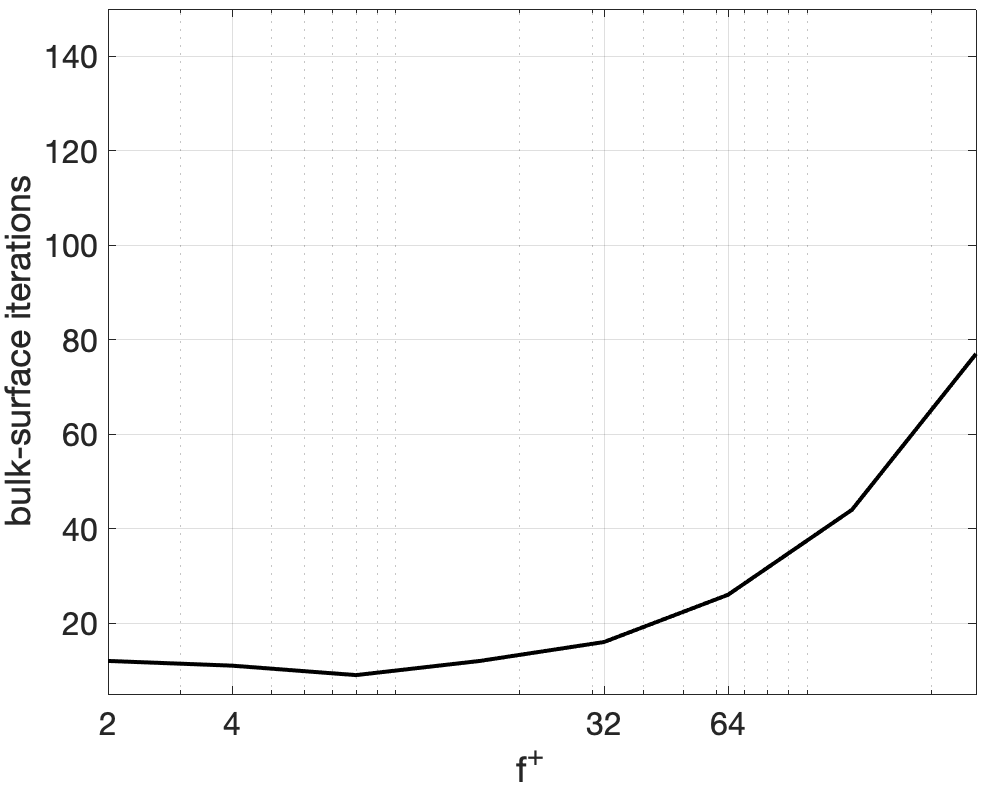}
  \caption{Sphere: (left) Bulk and surface FE errors
against the value of $f^+$. (right) Number of bulk-surface iterations of the partitioned method as $f^+$ varies.}
  \label{fig:fplus_sphere}
\end{figure}

\begin{figure}[hbt!]
\centering
\hskip -.3cm
\includegraphics[width=.56\textwidth]{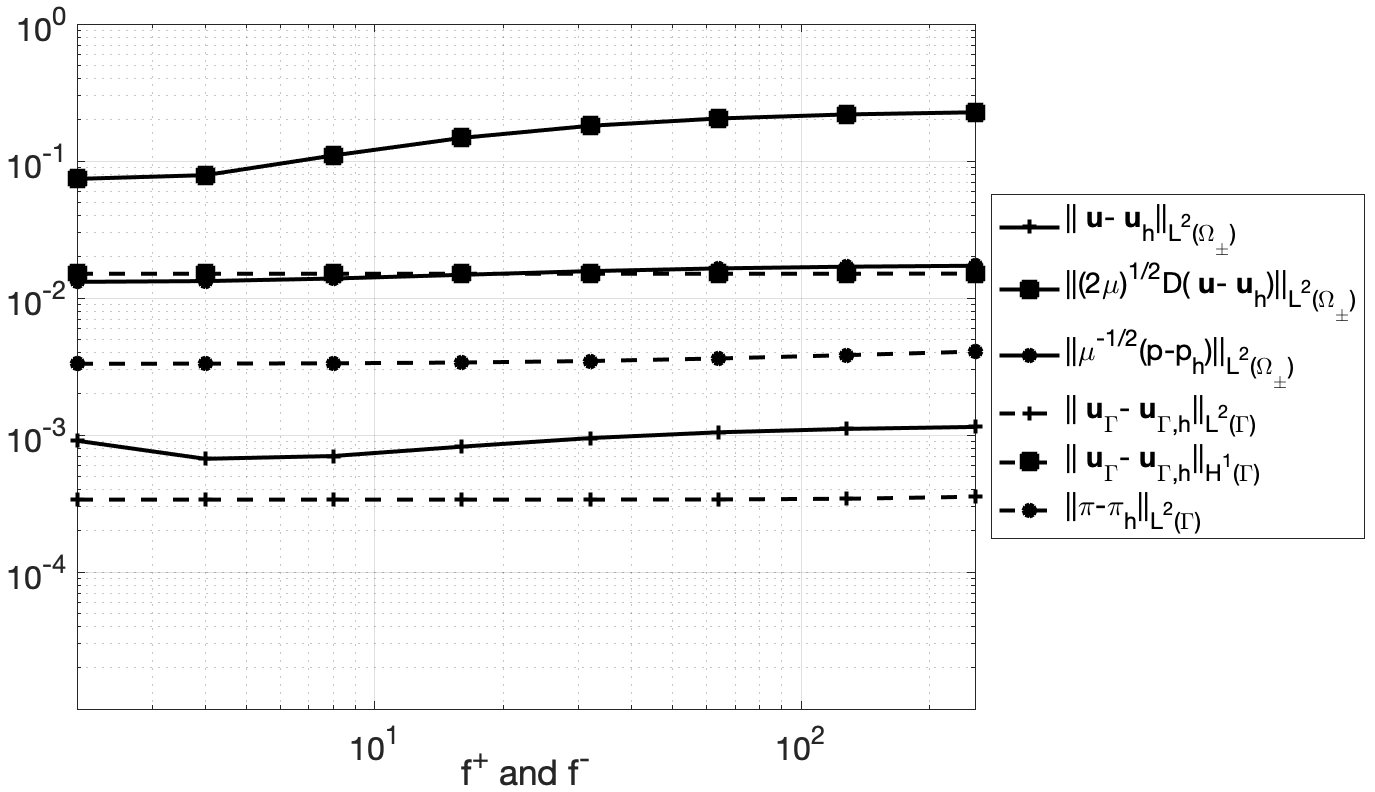}\quad
\includegraphics[width=.39\textwidth]{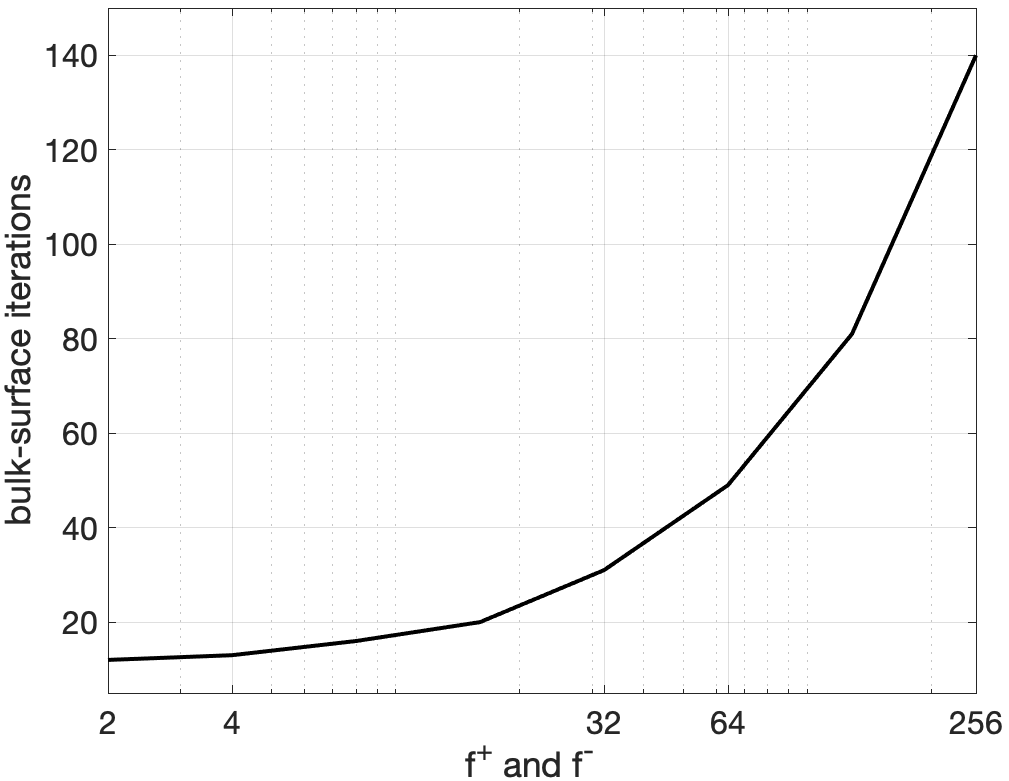}
  \caption{Sphere: (left) Bulk and surface FE errors
against the value of $f^+ = f^-$. (right) Number of bulk-surface iterations of the partitioned method as the value of $f^+$ and $f^-$
(with $f^+ = f^-$) varies.}
  \label{fig:f_sphere}
\end{figure}

{\bf Robustness with respect to the slip coefficients}. To check the sensitivity of the errors and partitioned
method to the value of the slip coefficients, we run two sets of experiments, both involving exact solution \eqref{ex_p_sphere}--\eqref{ex_surface_sphere}. In the first set we fix $f^+ = 2$ and vary $f^-$ from 1 to 256, while in the second set
we take $f^+ = f^-$ and let them both vary from 1 to 256 .
The viscosities are set as follows: $\mu^- = 1$,
$\mu^+ = 10$, and $\mu_\Gamma = 1$.
We consider again the mesh with mesh size $h = 0.125$.
Fig.~\ref{fig:fplus_sphere} (left) and \ref{fig:f_sphere} (left) show all the errors under consideration
against the value of the slip coefficient(s) for both sets of tests. The only error that shows a substantial
variation is the weighted $H^1$ error the bulk velocity, which increases as the slip coefficient(s) increase.
However, such error seems to reach a plateau in both cases.

\begin{figure}[htb]
    \centering
\includegraphics[width=.5\textwidth]{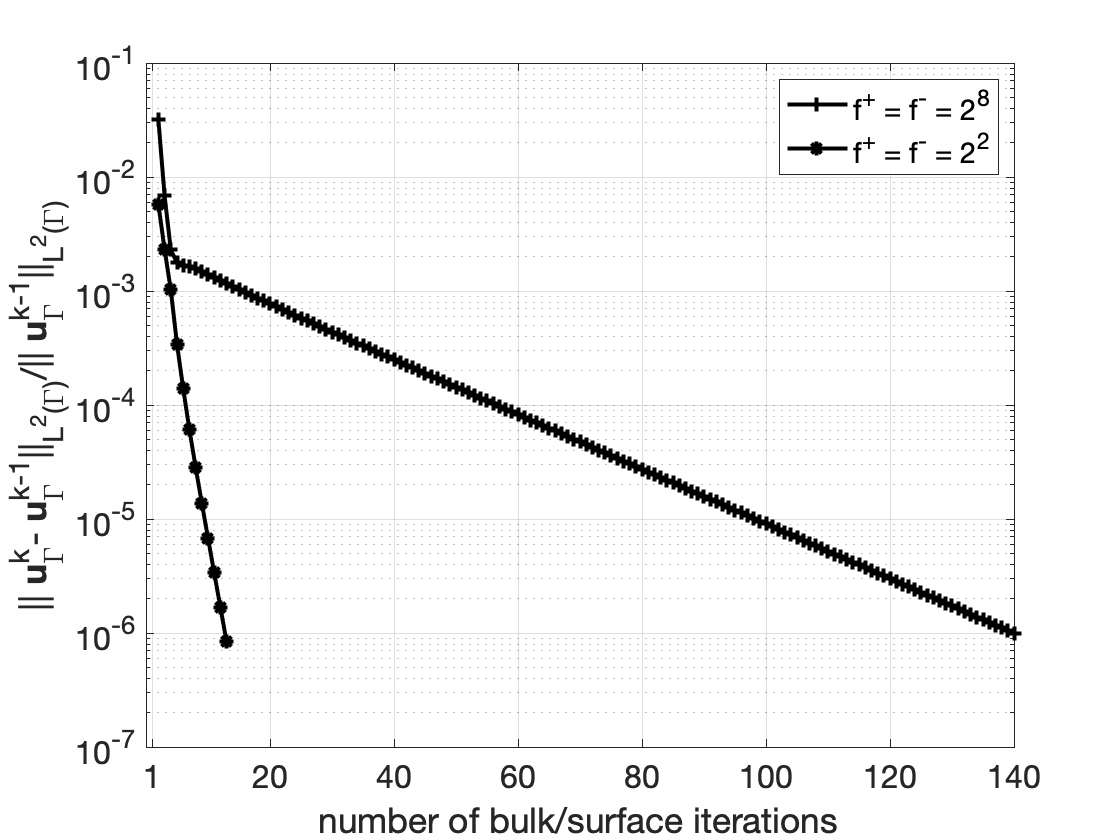}
    \caption{Sphere: relative difference of the surface velocity between subsequent iterations in $L^2$ norm until stopping criterion
    \eqref{eq:stop} is met.}
    \label{fig:iters}
\end{figure}

Fig.~\ref{fig:fplus_sphere} (right) and \ref{fig:f_sphere} (right) report the number of bulk-surface iterations to satisfy stopping
criterion \eqref{eq:stop} as the value of the coefficient(s) varies  for both sets of tests. In Fig.~\ref{fig:fplus_sphere} (right),
we see a rather sharp increase in the number of iterations as $f^-$ increases. This is even more true when both slip
coefficients are increased together, as we can see from Fig.~\ref{fig:f_sphere} (right).
Fig.~\ref{fig:iters} reports the relative difference of the surface velocity between subsequent iterations in $L^2$ norm
until stopping criterion \eqref{eq:stop} is met for $f^+ = f^- = 2^2$ and $f^+ = f^- = 2^8$. We see that such relative difference
decreases regularly for $f^+ = f^- = 2^2$, while for $f^+ = f^- = 2^8$ it decreases quickly for the first few iterations and then
it slows down. A heuristic explanation we have for this is that as the two friction coefficients increase interface
conditions \eqref{eq:cc2}--\eqref{eq:cc3} become close to Dirichlet conditions, making the surface flow more ``passive''.
Thus, separating the surface flow from the bulk flow as in the partitioned algorithm might not make much sense.

\subsection{Torus embedded in a cube}

The domain $\Omega$ is cube $[-2,2]^2$ and surface $\Gamma$
is a torus centered at $\bc = (c_1, c_2, c_3)$. Let
$(x,y)=(\tilde x-c_1,\tilde y-c_2, \tilde z - c_3)$, $(\tilde x,\tilde y, \tilde z)\in\Omega$.
We can characterize $\Gamma$ as the zero level set of function $\phi(\bx)= \sqrt{z^2+(\sqrt{x^2+y^2}-1)^2}-\frac{1}{2}$.
Finding an exact solution problem \eqref{eq:Stokes1}--\eqref{eq:scc4b},
\eqref{eq:cc2}, and \eqref{eq:cc3} with this more complicated surface is highly non-trivial. To simplify the task,
we relax interface conditions \eqref{eq:cc2}, \eqref{eq:cc3}, and \eqref{eq:scc4b} as follows:
\begin{align}
  \bP{\bsigma^\pm\bn}&=\pm f^\pm(\bP\bu^\pm-\bU )+g^\pm  &&\text{ on }\Gamma, \cl 
\left[ \bn^T \bsigma \bn \right]^-_+ &=\pi\kappa +g^n &&\text{ on }\Gamma. \el
\end{align}
where $g^+$, $g^-$, and $g^n$ are computed such that  exact solution given below satisfy these relaxed interface conditions.
The solution is given by
{\small
\begin{align}
p^-&=\left(\frac{1}{2}-\frac{2-4\sqrt{x^2+y^2}}{\sqrt{x^2+y^2}}\right)(x^3+x), ~ p^+=\frac{1}{2}(x^3+x), ~
\bu^-=\bu^+ =
\left[
\begin{array}{c}
x^2y\\
5-xy^2+z^2\\
-xy
\end{array}
\right], 
\label{ex_torus}
\end{align}
}
for the bulk and:
\begin{align}
\pi=x^3+x, \quad
\bU  =
\left[
\displaystyle{\frac{-zx}{\sqrt{x^2+y^2}}},~
\displaystyle{\frac{-zy}{\sqrt{x^2+y^2}}},~
\sqrt{x^2+y^2}-1
\right]^T, \label{ex_surface_torus}
\end{align}
for the surface.
The forcing terms $\bbf^-$ and $\bbf^+$ are found by plugging the solution \eqref{ex_torus}--\eqref{ex_surface_torus} in \eqref{eq:Stokes1}.
We impose a Dirichlet condition \eqref{eq:bcD} on the faces $x=2$, $y=-2$, $z=-2$, where function $\bg$
is found from $\bu^+$ in \eqref{ex_torus}. On the remaining part of the boundary, we impose
a Neumann condition \eqref{eq:bcN} where $\blf_n$ is found from $p^+$ and $\bu^+$ in \eqref{ex_torus}.

\begin{figure}[hbt!]
\centering
\hskip -.3cm
\includegraphics[width=.65\textwidth]{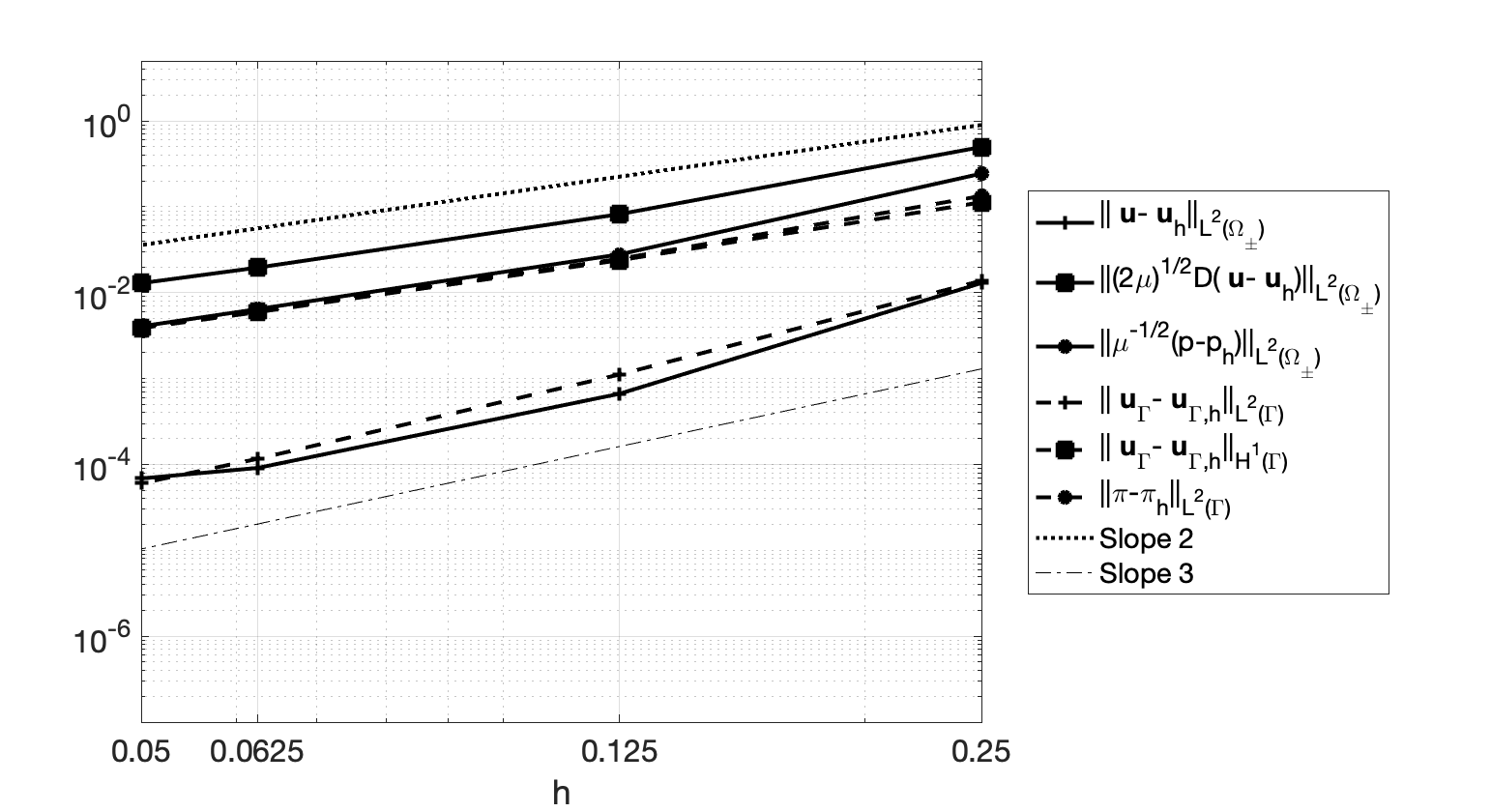}
  \caption{Torus: Bulk and surface FE errors against the mesh size $h$.}
  \label{fig:spatial_torus}
\end{figure}

\vskip .2cm
\noindent {\bf Spatial convergence.} Once again, we start by checking spatial accuracy. To this end, we consider
exact solution \eqref{ex_torus}--\eqref{ex_surface_torus} with $\bc = (0, 0, 0)$, viscosities $\mu^- = 1$,
$\mu^+ = 10$, $\mu_\Gamma = 1$, and friction coefficients $f^- = 2$ and $f^+ = 10$.
Just like in the case of the sphere, we consider structured meshes of tetrahedra that
feature a local one-level refinement near the corners of $\Omega$. The details
of the meshes under consideration are reported in Table \ref{tab:meshes_torus}.
Fig.~\ref{fig:spatial_torus} shows the $L^2$ error and weighted $H^1$ error
\eqref{eq:weighted_H1} for the bulk velocity, weighted $L^2$ error  \eqref{eq:weighted_L2}
for the bulk pressure,  $L^2$ and $H^1$ errors for the surface velocity and  $L^2$ error for the surface pressure
against the mesh size $h$. Also for this second convergence test, we
observe optimal convergence rates for all the norms.

\begin{table}[h]
\centering
\begin{tabular}{|c|c|c|c|c|}
\hline
$h$ & 0.25 & 0.125 & 0.0625 & 0.05  \\
\hline
\# bulk velocity DOFs & $1.6e5$ & $1.2e6$ & $8.5e6$ & $1.5e7$   \\
\hline
\# bulk pressure DOFs & $7.6e3$ & $5.4e4$ & $3.7e5$ & $6.7e5$  \\
\hline
\# surface velocity DOFs & $1.6e4$ & $6.0e4$ & $2.3e5$ & $3.4e5$ \\
\hline
\# surface pressure DOFs & $9.0e2$ & $3.4e3$ & $1.3e4$ & $2.0e4$ \\
\hline
\end{tabular}
\caption{Torus: DOFs for bulk and surface variables for all the meshes under consideration in the spatial convergence test.}
\label{tab:meshes_torus}
\end{table}

\begin{figure}[hbt!]
\centering
\hskip -.3cm
\includegraphics[width=.65\textwidth]{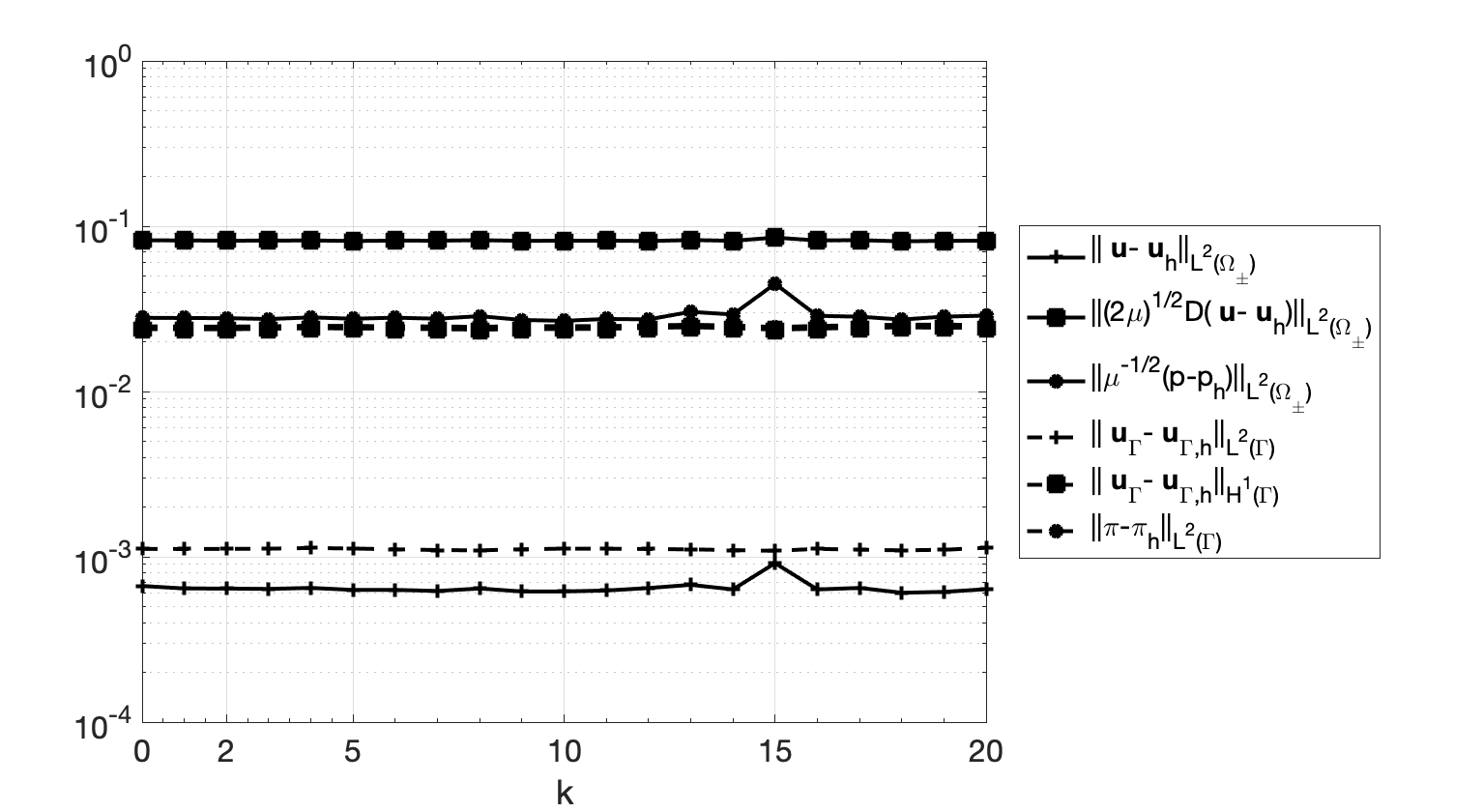}
  \caption{Torus: Bulk and surface FE errors
against the value of $k$ in \eqref{center}.}
  \label{fig:k_torus}
\end{figure}

{\bf Robustness with respect to the position of the interface}. We conclude our series of
numerical results with a set of simulations aimed at checking that our approach is not sensitive to
the position of the interface with respect to the background mesh. We vary
the center $\bc = (c_1, c_2, c_3)$ of the torus that represents $\Gamma$:
\begin{align}
c_1 = h \frac{k}{20} \sin\left(\frac{k\pi}{10}\right), ~c_2 = h \frac{k\sqrt{2}}{40} \cos\left(\frac{k\pi}{10}\right), ~c_3 = h \frac{k\sqrt{2}}{40} \cos\left(\frac{k\pi}{10}\right), \label{center}
\end{align}
where $h$ is the mesh size. The physical parameters are set like in the convergence test.
We consider the mesh in Table \ref{tab:meshes_torus} with $h = 0.125$.
Fig.~\ref{fig:k_torus} shows all the errors against the value of $k$ in \eqref{center}.
We see that all the errors are fairly insensitive to the position of $\Gamma$ with respect to the background mesh, indicating
robustness.

\section*{Acknowlegement} We are grateful to Dr. Christoph Lehrenfeld for providing us with an  ngsxfem implementation of isoparametric unfitted finite elements.  

\bibliographystyle{siam}
\bibliography{literatur}{}

\end{document}